\documentclass[10pt,draftcls,onecolumn]{IEEEtran}
\usepackage{graphicx}
\usepackage[ruled,vlined]{algorithm2e}
\usepackage{url}
 \usepackage[pdftex]{color}
\usepackage{amsmath}
\usepackage{amsfonts} 
\usepackage{array}
\usepackage{cancel}
\newtheorem{teor}{Theorem}[section]

\newtheorem{propo}{Proposition}[section]
\newtheorem{lemm}{Lemma}[section]

\usepackage{soul}
\usepackage{tabularx}
\usepackage{multirow}
\newcommand{\prop}{\begin{propo}}
\newcommand{\eprop}{\end{propo}}
\newcommand{\lem}{\begin{lemm}}
\newcommand{\elem}{\end{lemm}}
\newcommand{\teo}{\begin{teor}}
\newcommand{\eteo}{\end{teor}}
\newcommand{\ones}{\mathbf{1}_{\Sc}}
\newcommand{\tp}{^{\top}}

\newcommand{\beq}{\begin{equation}}
\newcommand{\eeq}{\end{equation}}
\newcommand{\bea}{\begin{eqnarray}}
\newcommand{\eea}{\end{eqnarray}}

\newcommand{\bsea}{\begin{subeqnarray}}
\newcommand{\esea}{\end{subeqnarray}}
\newcommand{\nn}{\nonumber}

\def\bmat{\left[ \begin{array}}
\def\emat{\end{array} \right]}
\DeclareMathOperator{\tr}{tr} 

\DeclareMathOperator*{\argmin}{arg\,min}
\definecolor{Royalblue}{cmyk}{1,0.30,0.2,0.2}
\newcommand{\al}[1]{\begin{align}#1\end{align}}

\newcommand{\Var}{\mathrm{Var}}

\newcommand{\Rs}{\mathbb{R}} 
\newcommand{\Es}{\mathbb{E}} 
\newcommand{\Nc}{\mathcal{N}} 
\newcommand{\Bc}{\mathcal{B}}

\newcommand{\Vc}{\mathcal{V}}
\newcommand{\Sc}{\mathcal{S}}
\newcommand{\Cc}{\mathcal{C}}

\newcommand{\col}[1]{\mathrm{col}(#1)}

\begin{document}

\title{Robust Distributed Kalman filtering with \\ Event-Triggered Communication}

\author{Davide~Ghion, Mattia~Zorzi\thanks{D. Ghion is with RMT s.r.l., Padova, Italy; e-mail:  {\tt\small d.ghion@ramete.com} .}\thanks{M. Zorzi is with the Department of Information
Engineering, University of Padova, Padova, Italy; e-mail: {\tt\small zorzimat@dei.unipd.it} .}}

\markboth{DRAFT}{}

\maketitle

\begin{abstract}
We consider the problem of distributed Kalman filtering for sensor networks in the case there are constraints in data transmission  and there is model uncertainty. More precisely, we propose two distributed filtering strategies with event-triggered communication where the state estimators are computed according to the least favorable model. The latter belongs to a ball  about the nominal model. We also show that both the methods are stable in the sense that the mean-square of the state estimation error is bounded in all the nodes.  
 \end{abstract}

\begin{IEEEkeywords}
Distributed robust Kalman filtering, sensor networks (SNs),
Event-triggered communication, model uncertainty.
\end{IEEEkeywords}

\section{Introduction}
Sensor networks are ubiquitous in many field, e.g. monitoring, security, data analysis and so on. The typical scenario is that  the sensors collect measurements and from them it is required to estimate some variables of interest, i.e. the state of a dynamical  model. This task is performed in a distributed fashion and it can be accomplished in different ways, see \cite{cattivelli2010,spanos2005,Saber_CDC2007,BATTISTELLI2014707,6576197,9091147,ryu2022consensus}; for instance, each sensor can update its estimate and then share the latter with its neighbors.

Nowadays the devices at each node of the sensor network are usually low-cost and battery-supplied. The latter feature guarantees that these devices are easily adaptable according to the surrounding environment. On the other hand, the data  transmission represents the most energy consuming node task. Therefore, it 
is fundamental to guarantee that the rate of the data   transmission across the network is below a certain threshold in order to avoid unexpected battery discharges. Data communication is typically scheduled by means of data-driven (or event-triggered) strategies, see  \cite{BATTISTELLI2012926,7047754,shi2016event},
\cite{Liu20152470,li2015event,yan2014distributed,shi2011time}. In this paper we focus our attention to the distributed scheme proposed in \cite{BATTISTELLI201875}: each node updates its estimate with the new measurement (if available); then, the latter is compared with the one just propagated in time through the state space model (i.e. the estimate that can be computed also by the neighbors of that node in the case there is no transmission). If the discrepancy between them is large, then the node sends out the updated estimate to its neighbors. Finally, each node performs a fusion of its estimate and the ones regarding its neighbors. 
The appealing property of this transmission rule is that it allows to obtain a distributed algorithm which enjoys nice stability properties (i.e., mean-square boundedness of the state estimation error in all the nodes) under minimal requirements.

In many situations the actual model is different from the nominal one (i.e. the one used in the estimation algorithm), \cite{RS_MPC_IET,9511139}. Indeed, the fact that the devices placed at each node are low-cost can lead to model uncertainty. In such a scenario, the performance of the estimator based on the nominal model will be poor. One possible way to address this issue is to consider the robust Kalman filter proposed in  \cite{ROBUST_STATE_SPACE_LEVY_NIKOUKHAH_2013}: the idea is to consider a dynamic minimax game where one player is the estimator which minimizes the state prediction error, while the other one is the hostile player which selects the least favorable model in a set of plausible models called ambiguity set. The latter is a ball which is formed by placing an upper bound on the discrepancy between the nominal state space model and the models inside. The radius of this ball is called tolerance and defines the magnitude of the uncertainty in the nominal model.

Distributed strategies which are robust to model uncertainty  have been already proposed in the literature, see for instance \cite{SHEN20101682,Luo_2008,RKDISTR,s20113244,9839566,9195416,caballero2020two}, however, to the best of the authors' knowledge, none of them consider the case with event-triggered communication.

In this paper we shall extend and formalize the ideas in the preliminary conference paper \cite{9837137}. More precisely, we present two new distributed Kalman filters with event-triggered communication under model uncertainty. The difference between the two approaches is the way with which the tolerance at each node collecting measurements is designed. In one approach the tolerance is the same for all those nodes, while in the other the tolerance depends on the local model corresponding to the node.  Our approaches represent a robust generalization of the distributed strategy proposed in \cite{BATTISTELLI201875}. We also show that, under reasonable assumptions,  both the methods are stable in the sense that the mean-square of the state estimation error is bounded in all the nodes. Finally, we also present a Monte Carlo study showing the effectiveness of the proposed approaches in the case of model uncertainty.

The outline of the paper is as follows. In Section \ref{sec:formulation} we formulate the distributed state estimation problem characterized by data transmission constraints and model uncertainty. In Section \ref{sec:algo} we introduce the robust distributed approach with event-triggered communication and uniform tolerance across the network. In Section \ref{sec_stab} we analyze the stability property of the proposed  method. In Section \ref{sec:local_algo} we propose the other robust distributed approach with event triggered communication. In Section \ref{sec:numerical} we consider the Monte Carlo study  showing the strength of our methods. Finally, in Section \ref{sec:conclusion} we draw the conclusions.

{\em Notation}. Given a matrix $A$, $A\tp$ denotes its transpose matrix. Given a symmetric matrix $Q$, then $Q>0$ ($Q\geq0$) means that $Q$ is positive (semi)definite; $\sigma_{min}(Q)$ and $\sigma_{max}(Q)$ denote the maximum and the minimum eigenvalue of $Q$, respectively. Given an index set $\Nc$ and a set of matrices $\{C^i, i\in \Nc\}$ having the same number of columns, then $\col{C_i,i\in\Nc}$  is the matrix obtained by stacking $C^i$'s. Given a set of matrices $\{R^i, i\in \Nc\}$, then $\mathrm{diag}(R^i,i\in\Nc)$ is the block diagonal matrix whose main blocks are $R^i$'s. Let $x$ be a random vector, then $\Es[x]$ and $\Var[x]$ denote its expectation and variance, respectively. Finally, $x  \sim\mathrm{N}(m,R)$ means that $x$ is Gaussian distributed with mean $m$ and covariance matrix $R$.


 \section{Problem formulation}\label{sec:formulation}
Consider a network of nodes described by the digraph $(\mathcal{N}, \mathcal{A}, \mathcal{S})$ where: $\mathcal{N} = \{1,...,N\}$ is the set of the nodes, $\mathcal{A} \subseteq \mathcal{N}\times\mathcal{N}$ is the set of edges and $\mathcal{S} \subseteq \mathcal{N}$ is the subset of sensors node. The latter are the only one that have the capabilities to perform  measurements, while the nodes in $\mathcal N\setminus \mathcal S$  are used to increase the connectivity of the network. If $(j,i)\in\mathcal A$, then it means that node $j$ can transmit data to node $i$; moreover, all the possible self-loops belong to $\mathcal A$. For each node $i\in\mathcal{N}$, the subset $\mathcal{N}_{i} := \{j:(j,i)\in\mathcal{A},\; j\neq i  \}\subseteq\mathcal{N}$ denotes the set of its in-neighbors, i.e. the set composed by the nodes that can send information to node $i$.  

We attach to this network the following nominal state space model
\begin{align}
      x_{t+1} &= Ax_{t}+Bw_{t}\label{eq:nominalStateEq}\\
    y^{i}_{t} & = C^{i}x_{t}+D^iv^{i}_{t}\;,\;\;\;\;  \;i\in\mathcal{S} \label{eq:nominalMeasEq}
\end{align}
where $x_{t}\in\mathbb R^n$ is the state, $y^i_t\in\mathbb R^{p_i}$ is the output at the sensor node $i$. Furthermore, matrix $B$ and $D^i$ are full row rank matrices, $w_{t}$ and $v^{i}_{t}$ are zero-mean normalized Gaussian white noises. The initial state $x_0$ is with mean $x_{0|-1}$ and with covariance matrix
$V_{0|-1}$. Finally, we assume that $w_{t}$, $v^{i}_{t}$'s and $x_0$ are independent. Notice that, the model (\ref{eq:nominalStateEq})-(\ref{eq:nominalMeasEq}) can be written as 
\begin{align}
      x_{t+1} &= Ax_{t}+Bw_{t}\label{eq:nominalStateEqglob}\\
    y_{t} & = Cx_{t}+Dv_{t} \label{eq:nominalMeasEqglob}
\end{align}
where $C:=\col{C^i, i\in \mathcal S}$,  $D:=\mathrm{diag}(D^i, i\in \mathcal S)$, $y_t:=\col{y_t^i, i\in \mathcal S}\in\mathbb R^p$ and $v_t:=\col{v_t^i, i\in \mathcal S}\in\mathbb R^p$.

The nominal model (\ref{eq:nominalStateEq})-(\ref{eq:nominalMeasEq}) over the time interval $t=1\ldots N$ can be equivalently described by the conditional probability densities $\phi_t(z_t|x_t)$, where $z_t=[\,x_{t+1}^\top\; y_t^\top\,]\tp$, $t=1\ldots N$, and $f(x_0)$. Taking the robust framework proposed in \cite{ROBUST_STATE_SPACE_LEVY_NIKOUKHAH_2013,STATETAU_2017,ZORZI2017133}, we assume that the actual model is described by the conditional probability densities $\tilde \phi_t(z_t|x_t)$, $t=1\ldots N$, and $f(x_0)$. Moreover, we assume that $\tilde \phi_t$ belongs to the ambiguity set 
$$\mathcal B_t=\{\,\tilde \phi_t \hbox{ s.t. } \tilde{\mathbb E}[\log(\tilde\phi_t/\phi_t) |Y_{t-1}]\leq b\,\},$$
\begin{align} &\tilde{\mathbb E}[\log(\tilde\phi_t/\phi_t) |Y_{t-1}]:=\nn\\
&\int_{\mathbb R^n}\int_{\mathbb R^{n+p}} \tilde \phi_t(z_t|x_t) \check f_t(x_t|Y_{t-1})\log\left( \frac{\tilde \phi_t(z_t|x_t) }{\phi_t(z_t|x_t) }\right)\mathrm d z_t \mathrm d x_t\nn,\end{align}
$Y_t:=\{y_0\ldots y_{t}\}$, $\check f_t$ is the actual probability density of $x_t$ given $Y_{t-1}$ and $b>0$ is called tolerance which accounts for model uncertainty. In plain words, $\mathcal B_t$ is a ball about the nominal density with with radius $b$. Then, a robust state estimator is the one solving the following minimax game:
\begin{align} x_{t+1|t}=\underset{g_t\in\mathcal G}{\argmin}\underset{\tilde \phi_t\in\mathcal B_t}{\max}\, \tilde{\mathbb E}[\| x_{t+1}-g_t(y_t)\|^2|Y_{t-1}]\label{central_pb}
\end{align}
where  $ x_{t+1|t}$ is the estimator of $x_{t+1}$ given $Y_{t}$; $\mathcal G$ is the set of estimators having finite second order moments for any $\tilde \phi_t\in \mathcal B_t$; \begin{align}
    \tilde{\mathbb{E}}&\left [ \|x_{t+1}-g_{t}(y_{t})\|^{2}|Y_{t-1} \right] :=\nn\\ & \int_{\mathbb R^{n}}\int_{\mathbb R^{n+p}}\|x_{t+1}-g_{t}(y_{t})\|^{2}\tilde{\phi}_{t}(z_{t}|x_{t})\check{f}_{t}(x_{t}|Y_{t-1})dz_{t}dx_{t}\nn
\end{align}
and it is assumed that $\check{f}_{t}(x_{t}|Y_{t-1})\sim\mathrm{N}(\hat{x}_{t|t-1}, V_{t|t-1})$. The basic idea behind this paradigm is that when we are looking for an estimator that minimizes properly the selected loss function, a hostile player called ``nature'' conspires to select the worst possible model in the ambiguity set $\mathcal{B}_{t}$. In \cite{ROBUST_STATE_SPACE_LEVY_NIKOUKHAH_2013} it has been shown that the (centralized) robust estimator solution to (\ref{central_pb}) admits a Kalman-like structure. It is not difficult to show that such a filter can be written in the information form as follows. Let $P_{t|t-1}$ and $V_{t|t-1}$ denote the pseudo-nominal and least favorable, respectively, covariance matrix of the prediction error at time $t$; let $x_{t|t}$ denote the estimator of $x_t$ given $Y_t$ and $P_{t|t}$ denotes the covariance matrices of the corresponding estimation error at time $t$.  We define the corresponding matrices in the information form as $\Omega_{t|t-1} = P_{t|t-1}^{-1}$, $\Psi_{t|t-1} = V_{t|t-1}^{-1}$, $\Omega_{t|t} = P_{t|t}^{-1}$ and the information states as
\begin{equation*}
q_{t|t-1}=\Psi_{t|t-1}{x}_{t|t-1},\quad q_{t|t}=\Omega_{t|t}{x}_{t|t}.
\end{equation*} Then, it is not difficult to see that the robust estimator obeys to
\begin{align}
      {\small \text{Correction step:}}&
    \begin{cases}
        \Omega_{t|t} = \Psi_{t|t-1}+ C\tp R^{-1}C\nn\\
        q_{t|t} = q_{t|t-1} + C\tp R^{-1}y_{t} ,      
    \end{cases}\\
       {\small\text{Prediction step:}}&
    \begin{cases} 
    \Omega_{t+1|t} =\\ \hspace{0.2cm}     Q^{-1}-Q^{-1}A(A\tp Q^{-1}A+\Omega_{t|t})A\tp Q^{-1}\nn\\
                \text{Find } \theta_{t}>0 \text{ s.t. }\gamma(\Omega_{t+1|t}, \theta_{t})=b\\
  \Psi_{t+1|t} = \Omega_{t+1|t}-\theta_{t}\mathbf{I}_{n}\\
        q_{t+1|t} = \Psi_{t+1|t}A\Omega^{-1}_{t|t}q_{t|t},    
    \end{cases}
\end{align}
 where $R:=DD\tp$, $Q:=BB\tp$ are positive definite matrices and
 \al{ \gamma&(\Omega,\theta) \nn\\ &:= \frac{1}{2}\left\{ tr[(\mathbf{I}_{n}-\theta \Omega^{-1})^{-1}-\mathbf{I}_{n}] + \log\det(\mathbf{I}_{n}-\theta \Omega^{-1})\right \}.\nn }The parameter $\theta_{t}>0$ is called risk sensitivity parameter. It is worth noting that given $\Omega>0$ and $b>0$, the equation $\gamma(\Omega,\theta)=b$ always admits a unique solution $\theta>0$ such that $\Omega-\theta \mathbf{I}_{n}>0$. Furthermore, in the special case where $b =0$, i.e. there is no model uncertainty, we have that $\gamma(\Omega,\theta)=0$ implies that $\theta_{t}=0$ and thus the above equations degenerates in the usual Kalman equations in the information form.
 
In what follows, we face the problem to solve the minimax game in (\ref{central_pb}) where the minimizer works in a distributed way and under data transmission constraints. More precisely, each node $i\in\mathcal N$ must estimate the state $x_t$, with $t\in\mathbb{Z}_{+}=\{1, ..., N\}$,   taking into account that: 
i) the actual model does not coincide with the nominal one (\ref{eq:nominalStateEqglob})-(\ref{eq:nominalMeasEqglob}); ii) each node $i$ can selectively transmits only the most relevant data without compromising the stability properties.
 

\section{Robust event-triggered strategy}\label{sec:algo}
Before to introduce our robust distributed estimation paradigm we consider the following quite simple scenario. We assume that $\mathcal N=\mathcal S$, i.e. all the nodes are sensor nodes, and $\mathcal A=\{\,(j,j),\; j\in\Nc\,\}$, i.e. the nodes do not communicate. In the presence of model uncertainty, at node $i\in\Sc$ we can consider the robust Kalman filter in the information form based on the local model
 \begin{align}
      x_{t+1} &= Ax_{t}+Bw_{t}\nn\\
    y^{i}_{t} & = C^{i}x_{t}+D^iv^{i}_{t}\;\nn ;   
\end{align}
thus, we obtain the following algorithm:
\begin{align}
     {\small \text{Correction step:}}&
    \begin{cases}\label{eq:RKFfilteringloc}
        \Omega_{t|t}^i = \Psi_{t|t-1}^i+ (C^i)\tp(R^i)^{-1}C^i\\
        q_{t|t}^i = q_{t|t-1}^i + (C^i)\tp(R^i)^{-1}y_{t}^i       
    \end{cases}\\
    {\small \text{Prediction step:}}&
    \begin{cases}\label{eq:RKFpredictionloc}
    \Omega_{t+1|t}^i =\\ \hspace{0.2cm}     Q^{-1}-Q^{-1}A(A\tp Q^{-1}A+\Omega_{t|t}^i)A\tp Q^{-1}\\
                \text{Find } \theta_{t}^i>0 \text{ s.t. }\gamma(\Omega_{t+1|t}^i, \theta_{t}^i)=b\\
  \Psi_{t+1|t}^i = \Omega_{t+1|t}^i-\theta_{t}^i\mathbf{I}_{n}\\
        q_{t+1|t}^i= \Psi_{t+1|t}^iA(\Omega^{i}_{t|t})^{-1}q_{t|t}^i
    \end{cases}
\end{align}
 where $R^i:=D^i(D^i)\tp$. It is worth noting that each node has its own risk sensitivity parameter $\theta_{t}^i$. 
  
Next, we consider the scenario in which: i) $ \mathcal S\subseteq \mathcal N$; ii) each node $i$ can transmit its local estimate ${q}^{i}_{t|t}$ and information matrix $\Omega^{i}_{t|t}$ to all its out-neighbors $\overline{\Nc}_i=\{\,j\, : \, (i,j)\in\mathcal A, \, i\neq j\,\}$ if necessary.  In plain words, each node $i$ can decide at any time step whether to transmit or not its data, i.e. $(q^i_{t|t},\Omega^{i}_{t|t})$,  without compromising the stability properties of the algorithm. More precisely, we require that $\Es[\|x_t-x^i_{t|t}\|^2]$, where the expectation operator is under the global least favorable model solution to (\ref{central_pb}), does not diverge for any $i\in\mathcal N$ as $t$ approaches infinity.

The estimation paradigm that we now present  
is a robust generalization of the distributed state estimation algorithm with event-triggered communication proposed in  \cite{BATTISTELLI201875} and it is composed by four steps described below.

\textbf{Correction.} At time $t$,  the predicted pair $(q^{i}_{t|t-1}, \Psi^{i}_{t|t-1})$ is available at node $i\in \mathcal N$. If $i\in\mathcal S$, i.e. it is a sensor node, then also the measurement $y_t^i$ is available and thus the correction step coincides with (\ref{eq:RKFfilteringloc}). If $i\notin\mathcal S$,
no measurement is available at the node, then we can only   propagate the prediction couple. Therefore, the so called information pair is obtained as
\begin{equation} 
    (q^{i}_{t|t},\Omega^{i}_{t|t}) =
    \begin{cases}
        \text{use equations (\ref{eq:RKFfilteringloc}),} & \text{if $i\in\mathcal{S}$}\nn\\
        (q^{i}_{t|t-1}, \Psi^{i}_{t|t-1}), & \text{if $i\in\mathcal{N \setminus 
        S}$}.
    \end{cases}
\end{equation}

\textbf{Information exchange.} Each node $i\in\mathcal N$ sends its information couple $(q^{i}_{t|t}, \Omega^{i}_{t|t})$ to all its out-neighbors according to the binary variable  $c^{i}_{t}$: \begin{itemize} 
\item if $c^{i}_{t}=1$, then node $i $ transmits the information couple to all its out-neighbors at time $t$;
 \item  if $c^{i}_{t}=0$, then node $i $ does not transmit the  information couple to all its out-neighbors at time $t$.
 \end{itemize}
 It remains to define the binary variable $c^i_t$. Let $n^{i}_{t}\in\mathbb N$ be the number of time instants elapsed from the most recent transmission of node $i$, i.e.  the most recently transmitted data is $(q^{i}_{t-n^{i}_{r}|r-n^{i}_{t}}, \Omega^{i}_{t-n^{i}_{t}|t-n^{i}_{t}})$. Then, all the out-neighbors of node $i$ propagate $(q^{i}_{t-n^{i}_{r}|t-n^{i}_{t}}, \Omega^{i}_{t-n^{i}_{t}|t-n^{i}_{t}})$ in time through a prediction step which takes into account the fact that the actual model does not coincide with the nominal one (see (\ref{eq:estPredictionStep}) in the prediction step below). Let  $(\bar{q}^{i}_{t}, \bar{\Psi}^{i}_{t})$ denote this propagated pair at time $t$. Then, the transmission rule $c_t^i$ computed at node $i$ is defined as in \cite{BATTISTELLI201875}:  
 \begin{equation}\label{eq:eventTriggered}
    c^{i}_{t} = 
    \begin{cases}
        0, & \quad \text{if } \| {x}^{i}_{t|t}-\bar{x}^{i}_{t}\|^{2}_{\Omega^{i}_{t|t}}\leq\alpha\\ & \text{ and } \frac{1}{1+\beta}\Omega^{i}_{t|t}\leq \bar{\Psi}^{i}_{t} \leq (1+\delta)\Omega^{i}_{t|t}\\
        1, & \quad \text{otherwise}
    \end{cases}
\end{equation}
where $\bar{x}^{i}_{t} =  (\bar{\Psi}^{i}_{t})^{-1}\bar{q}^{i}_{t}$ represents the robust state prediction based on the propagation of the most recent transmitted pair $(q^{i}_{t-n^{i}_{r}|t-n^{i}_{t}}, \Omega^{i}_{t-n^{i}_{t}|t-n^{i}_{t}})$; $x^{i}_{t|t} = ({\Omega}^{i}_{t|t})^{-1}q^{i}_{t|t}$ is the state estimate at node $i$. In plain words, the transmission rule in (\ref{eq:eventTriggered})
checks the discrepancy between $(q^{i}_{t|t}, \Omega^{i}_{t|t})$ and $(\bar{q}^{i}_{t}, \bar{\Psi}^{i}_{t})$. If the latter is large,
then it means that the out-neighbors own a prediction corresponding to node $i$ which is bad and thus node $i$ must transmit the data.   
The positive scalars $\alpha$, $\beta$ and $\delta$ can be tuned by the user in order to reach a desired behavior in terms of transmission rate and performance. More precisely, $\alpha$ tunes the bound on the discrepancy between $x^i_{t|t}$ and $\bar x^i_t$, while $\beta$ and $\delta$ tunes the allowed mismatch between the covariance matrices ${\Omega}^i_{t|t}$ and $\bar{\Psi}^i_t$. In   \cite{BATTISTELLI16} it has been shown that   
the transmission strategy in (\ref{eq:eventTriggered}) guarantees the following upper bound. If we model the propagated and the information pairs as $\mathrm N(\bar{q}^{i}_{t},  \bar{\Psi}^{i}_{t})$ and $\mathrm N (q^{i}_{t|t}, \Omega^{i}_{t|t})$,  respectively,  then Condition (\ref{eq:eventTriggered}) guarantees that  \begin{equation*} 
    D_{KL}(\mathbb N(q^{i}_{t|t}, \Omega^{i}_{t|t}),\mathbb N(\bar{q}^{i}_{t}, \bar{\Psi}^{i}_{t})) \leq \frac{1}{2}[\alpha+\beta n+n\log(1+\delta)]
\end{equation*}
where $n$ is the state dimension and $D_{KL}$ denotes the Kullback-Leibler divergence.

\textbf{Information fusion.} In this step, any node merges its information with the ones corresponding to its in-neighbors. Let $\Pi\in\mathbb R^{N\times N}$ denote the consensus matrix whose element
in position $(i,j)$ is defined as:
\begin{equation*}
    \pi_{i,j} = \begin{cases}
          (d_{i}+1) ^{-1}, & \text{ if }(j,i)\in\mathcal{A}\\
          0, & \text{otherwise}
    \end{cases}
\end{equation*} where $d_i$ denotes the degree of node $i$; in this way we have that 
$\pi_{i,j}$ with $j\in\mathcal N_i$ represents the coefficients of a convex combination. 
 Then, the fusion step is performed through the following convex combination of the  pairs:
 \begin{equation}\label{eq:qFusedwithTilde}
    q^{i,F}_{t|t} = \pi_{i,i}q^{i}_{t|t}+\sum_{j\in\mathcal{N}_{i}}\pi_{i,j}\left [c^{j}_{t}q^{j}_{t|t}+(1-c^{j}_{t})\tilde{q}^{j}_{t}\right ]
\end{equation}
\begin{equation}\label{eq:OmegaFusedwithTilde}
    \Omega^{i,F}_{t|t} = \pi_{i,i}\Omega^{i}_{t|t}+\sum_{j\in\mathcal{N}_{i}}\pi_{i,j}\left [c^{j}_{t}\Omega^{j}_{t|t}+(1-c^{j}_{t})\tilde{\Omega}^{j}_{t}\right ] 
\end{equation}
where 
\begin{equation*}
    \tilde{q}^{j}_{t} = \frac{1}{1+\delta}\bar{q}^{j}_{t}\;\;\;,\;\;\;\tilde{\Omega}^{j}_{t} = \frac{1}{1+\delta} \bar{\Psi}^{j}_{t}.
\end{equation*}
In view of (\ref{eq:qFusedwithTilde})-(\ref{eq:OmegaFusedwithTilde}), we can see that in the fusion step we consider $(q^j_{t|t},\Omega^j_{t|t})$, if node $j$ transmitted its information pair at time $t$. If node $j$ does not transmit, then the aforementioned pair is not available at node $i$. To account for this lack, given that at each iteration the nodes can calculate the pair $(\bar{q}^{i}_{t}, \bar{\Psi}^{i}_{t})$, which is certainly less informative than $(q^{i}_{t|t}, \Omega^{i}_{t|t})$, we consider $(\bar{q}^{i}_{t}, \bar{\Psi}^{i}_{t})$ in  (\ref{eq:qFusedwithTilde})-(\ref{eq:OmegaFusedwithTilde}) shrunk by   the factor $(1+\delta)^{-1}$ in order to decrease its importance in the fusion step, see \cite{BATTISTELLI201875} for more details.

\textbf{Prediction.} Once each node $i\in\mathcal{N}$ has computed the fused information couple $(q^{i,F}_{t|t} \Omega^{i,F}_{t|t})$, the latter is propagated in  time with the robust prediction step in (\ref{eq:RKFpredictionloc}) where $(q^i_{t|t},\Omega_{t|t}^{i})$ is now replaced by $(q^{i,F}_{t|t},\Omega_{t|t}^{i,F})$:
\begin{align}\label{eq:infoFormPred} 
\left.\begin{array}{ll}
    \Omega_{t+1|t}^i &=Q^{-1}-Q^{-1}A(A\tp Q^{-1}A+\Omega_{t|t}^{i,F})A\tp Q^{-1}\\
                \text{Find } &\theta_{t}^i>0 \text{ s.t. }\gamma(\Omega_{t+1|t}^i, \theta_{t}^i)=b\\
  \Psi_{t+1|t}^i &= \Omega_{t+1|t}^i-\theta_{t}^i\mathbf{I}_{n}\\
        q_{t+1|t}^i&= \Psi_{t+1|t}^iA(\Omega^{i,F}_{t|t})^{-1}q_{t|t}^{i,F}.
\end{array}\right. 
\end{align}
In this step we compute the propagated pair
$(\bar{q}^{i}_{t+1},  \bar{\Psi}^{i}_{t+1} )$, i.e. the one used in the case node $i$ does not transmit its information pair. 
Notice this operation is performed by both node $i$ and its out-neighbors and it can be summarized as follows.  At nodes $i\,\cup\, \mathcal N_i $  we have the pair 
$(\breve{q}^{i}_{t},\breve{\Omega}^{i}_{t})$ defined as
\begin{align}
\begin{cases}\label{def_breve}
\breve{q}^{i}_{t} &= c^{i}_{t}q^{i}_{t|t}+(1-c^{i}_{t})\bar{q}^{i}_{t}\\ 
\breve{\Omega}^{i}_{t} &= c^{i}_{t}\Omega^{i}_{t|t}+(1-c^{i}_{t})\bar{\Psi}^{i}_{t};
\end{cases}
\end{align} then, it is propagated in  time with  the robust prediction step in (\ref{eq:RKFpredictionloc}) where $(q^i_{t|t},\Omega_{t|t}^{i})$ is now replaced by $(\breve q^{i}_{t},\breve \Omega_{t}^{i})$: 
\begin{align}
\label{eq:estPredictionStep}
\left.\begin{array}{ll}
   \bar{ \Omega}_{t+1}^i &=Q^{-1}-Q^{-1}A(A\tp Q^{-1}A+\breve\Omega_{t}^{i})A\tp Q^{-1}\\
                \text{Find } &\bar\theta_{t}^i>0 \text{ s.t. }\gamma(\bar \Omega_{t+1}^i, \bar\theta_{t}^i)=b\\
 \bar \Psi_{t+1}^i &= \bar\Omega_{t+1}^i-\bar\theta_{t}^i\mathbf{I}_{n}\\
        \bar q_{t+1}^i&= \bar \Psi_{t+1}^iA(\breve\Omega^{i}_{t})^{-1}\breve q_{t}^{i}.
\end{array}\right. 
\end{align}

The procedure is summarized in Algorithm \ref{alg:DKFeventTriggered}. It is worth noting that each node $i$ is characterized by two risk sensitivity parameters, i.e. $\theta^i_t$ and $\bar \theta^i_t$. In the case that $b=0$, i.e. there is no model uncertainty, in Algorithm \ref{alg:DKFeventTriggered} we have: $\theta^i_t=0$, $\bar \theta^i_t=0$ and thus $\Psi_{t}^i = \Omega_{t}^i$, $ \bar \Psi_{t}^i = \bar\Omega_{t}^i$, i.e. we recover the distributed Kalman algorithm with event-triggered communication proposed in \cite{BATTISTELLI201875}.
\begin{algorithm}
    \caption{\textbf{RDKF} with event-triggered communication}\label{alg:DKFeventTriggered}
    \textbf{Initialization:} Set $(q^{i}_{0|-1}, \Psi^{i}_{0|-1})$ for any $i\in\mathcal N$
    \newline
    For each $t=0,1,\ldots$
    \newline
    For each node $i\in\mathcal{N}$
    \begin{description}
        \item \textbf{Correction:} 
        \begin{equation*}
    (q^{i}_{t|t},\Omega^{i}_{t|t}) =
    \begin{cases}
        \text{use   (\ref{eq:RKFfilteringloc}),} & \text{if $i\in\mathcal{S}$}\\
        (q^{i}_{t|t-1}, \Psi^{i}_{t|t-1}), & \text{if $i\in\mathcal{N \setminus 
        S}$}
    \end{cases}
\end{equation*}
        \item \textbf{Information exchange:} 
        \begin{description}
        \item[-] if $t=0$ set $c^i_t=1$, otherwise determine $c_t^i$ according to (\ref{eq:eventTriggered})
        \item[-] if $c^i_t=1$ transmit $(q^i_{t|t},\Omega^i_{t|t})$ to the out-neighbors 
        \item[-] receive $(q^j_{t|t},\Omega^j_{t|t})$ from all the in-neighbors $j\in\mathcal N_i$ for which $c_t^j=1$
        \end{description}
\item \textbf{Information fusion:} 
\begin{equation*}
\hspace*{-0cm}    \tilde{q}^{j}_{t} = \frac{1}{1+\delta}\bar{q}^{j}_{t}\;\;\;,\;\;\;\tilde{\Omega}^{j}_{t} = \frac{1}{1+\delta} \bar{\Psi}^{j}_{t} \text{ with } j\in\mathcal N_i
\end{equation*}
 \begin{align*}
    q^{i,F}_{t|t} &= \pi_{i,i}q^{i}_{t|t}+\sum_{j\in\mathcal{N}_{i}}\pi_{i,j}\left [c^{j}_{t}q^{j}_{t|t}+(1-c^{j}_{t})\tilde{q}^{j}_{t}\right ]\\
    \Omega^{i,F}_{t|t} &= \pi_{i,i}\Omega^{i}_{t|t}+\sum_{j\in\mathcal{N}_{i}}\pi_{i,j}\left [c^{j}_{t}\Omega^{j}_{t|t}+(1-c^{j}_{t})\tilde{\Omega}^{j}_{t}\right ] 
\end{align*}
\item \textbf{Prediction step:} 
 \begin{description}
 \item[-] Compute $(q^i_{t+1|t},\Psi^i_{t+1|t})$ using (\ref{eq:infoFormPred})
 \item[-] Compute $(\bar q^i_{t+1|t},\bar \Psi^i_{t+1})$ using (\ref{def_breve}) and (\ref{eq:estPredictionStep})
 \end{description}
     \end{description}
\end{algorithm}
\section{Stability analysis} \label{sec_stab}
Recall that the nominal (global) model (\ref{eq:nominalStateEqglob})-(\ref{eq:nominalMeasEqglob}) is different from the actual one. In this section we analyze the stability properties of Algorithm \ref{alg:DKFeventTriggered} under the least favorable model which is given by the maximizer of the minimax problem in (\ref{central_pb}), that is, the centralized problem. In doing that, we need the following assumptions:
\begin{description}
\item[\textbf{A1}.] The tolerance $b$ defining the ambiguity set in (\ref{central_pb}) is taken sufficiently small;
\item[\textbf{A2}.] The transition matrix A is invertible;
\item[\textbf{A3}.] The system is collectively observable, i.e., the pair $(A,C)$ is observable;
\item[\textbf{A4}.] The network is strongly connected, that is, there exists a directed path between any pair $i,j\in\Nc$.
\end{description}
It is worth noting that assumptions \textbf{A3}-\textbf{A4} are the same made in \cite{BATTISTELLI201875}, in particular Assumption \textbf{A2} automatically holds in sampled-data systems where matrix $A$ is obtained by discretizing the corresponding continuous-time matrix. Finally, Assumption \textbf{A1} is necessary to guarantee that the noise processes characterizing the global least favorable model have uniformly bounded  variance, see the next proposition.

\prop\label{propo_leastfav}
Under assumptions \textbf{A1} and \textbf{A3}, the least favorable model solution to (\ref{central_pb}) takes the form
\begin{align}
      x_{t+1} &= Ax_{t}+B\tilde{w}_{t}\label{eq:nominalStateEqgloblf}\\
    y_{t} & = Cx_{t}+D\tilde{v}_{t} \label{eq:nominalMeasEqgloblf}
\end{align}
where $\tilde w_t$ and $\tilde v_t:=\col{\tilde v_t^i,i\in\Sc}$ are zero-mean Gaussian colored noises such that  
\begin{align*}
\underline{\rho} I_{n+N}\leq \Var\left(\left[\begin{array}{c}
\tilde w_t\\ \tilde v_{t+1}\end{array}\right]\right)\leq \overline{\rho} I_{n+N}
\end{align*}
with $\overline{\rho} \geq \underline{\rho} \geq 0$. Moreover, $\tilde w_t$ and $\tilde v_t$ are correlated. 
\eprop 
\begin{proof} The least favorable model (\ref{eq:nominalStateEqgloblf})-(\ref{eq:nominalMeasEqgloblf}) has been characterized in \cite[Section V]{ROBUST_STATE_SPACE_LEVY_NIKOUKHAH_2013} where 
\al{&\tilde w_t:=H_{1,t}\varepsilon_t+L_{1,t}\epsilon_t, \quad \tilde v_t:=H_{2,t} \varepsilon_t+L_{2,t} \epsilon_t;\nn}
$\epsilon_t$ is zero-mean normalized Gaussian white noise,  \al{\varepsilon_{t+1}=[A+ BH_{1,t}-&G_t (C+DH_{2,t})]\varepsilon_t\nn\\ &+(BL_{1,t}-G_tDL_{2,t})\epsilon_t;\nn}
the definition of the matrices $H_{1,t}$, $H_{2,t}$, $L_{1,t}$, $L_{2,t}$ and $G_t$ can be found in \cite{ROBUST_STATE_SPACE_LEVY_NIKOUKHAH_2013}. In \cite{zorzi2018lf} it was shown that, under Assumptions \textbf{A1} and \textbf{A3}, $H_{1,t}\rightarrow H_1$, $H_{2,t}\rightarrow H_2$, $L_{1,t}\rightarrow L_1$, $L_{2,t}\rightarrow L_2$ and  $G_t\rightarrow G$ as $t\rightarrow \infty$. Moreover,  $\check A:=A+ BH_1- G (C+DH_2)$ is Schur stable. Accordingly $\Es[\varepsilon_t \varepsilon_t\tp]\rightarrow E_0$, with $E_0\geq0$, and $\Es[\varepsilon_{t+1} \varepsilon_t\tp]\rightarrow \check E_1=\check AE_0$ as $t\rightarrow \infty$. Then,
\al{\Var &\left(\left[\begin{array}{c}
\tilde w_t\\ \tilde v_{t+1}\end{array}\right]\right)
&\rightarrow HEH\tp+LL\tp+HM L\tp+LM\tp H\tp\nn}
where
 \al{ & H:=\left[\begin{array}{cc}
H_{1}& 0\\ 0& H_{2}\end{array}\right],\quad 
E:=\left[\begin{array}{cc}
E_{0}& E_1\tp\\ E_1& E_{0}\end{array}\right],\nn\\
& L:=\left[\begin{array}{cc}
L_{1}& 0\\ 0& L_{2}\end{array}\right],\quad
M:=\left[\begin{array}{cc}
0& 0\\ BL_{1}-GDL_{2}& 0\end{array}\right]. \nn}   
Hence, we have
\al{\underline{\rho}&=\sigma_{min}(HEH\tp+LL\tp+HM L\tp+LM\tp H\tp),\nn\\  \overline{\rho}&=\sigma_{max}(HEH\tp+LL\tp+HM L\tp+LM\tp H\tp).\nn}
\end{proof}
\smallskip
Let $e_t^i=x_t-x^i_{t|t}$ denote the estimation error at node $i$ and $e_t=\col{e_t^i,i\in\Nc}$ the collective estimation error. In order to prove the stability of the estimation error, we consider the Lyapunov function
\al{ \Vc_t(e_t)=\sum_{i\in\Nc} p_i\|e_t^i\|^2_{\Omega_{t|t}^i}} where $p_i$'s, strictly positive, are the components of a vector $p$ satisfying the condition $p\tp=p\tp\Pi$. Notice that, the existence of such a vector is guaranteed by the Perron-Frobenious theorem  since, by Assumption \textbf{A4}, $\Pi$ is a primitive matrix. The next proposition shows that $\Vc_t(e_t)$ is a well-defined Lyapunov function because $\Omega_{t|t}^i$'s are uniformly bounded.
\smallskip
\prop \label{prop_bound_Omega}Assume that \textbf{A2}-\textbf{A4} hold. Consider the sequence $\Omega_{t|t}^i$, $i\in\Nc$, generated by Algorithm \ref{alg:DKFeventTriggered} with $\Psi_{0|-1}^i>0$. Then, there exist three positive constants $\underline{\omega}$, $\overline{\omega}$ and $\underline{\omega}_F$ such that $\underline{\omega}I\leq  \Omega_{t|t}^i\leq\overline{\omega}I$ and $\underline{\omega}I\leq  \Omega_{t|t}^{i,F}\leq\overline{\omega}_FI$.
\eprop 
\begin{proof} We start by showing the upper bounds. By (\ref{eq:RKFpredictionloc}) we have that $\Omega_{t+1|t}^i\leq Q^{-1}$ and thus $\Psi_{t+1|t}^i\leq \Omega_{t+1|t}^i\leq Q^{-1}$ which is a uniform upper bound. By the correction step in Algorithm \ref{alg:DKFeventTriggered} we have that $\Omega^{i}_{t|t}\leq Q^{-1}+(C^i)\tp(R^i)^{-1}C^i\leq \overline{\omega}I$
where 
\al{\overline\omega:=\underset{i\in\Nc}{\max}\,\sigma_{max}(Q^{-1}+(C^i)\tp(R^i)^{-1}C^i).\nn} It is also worth noting that \al{\label{upper_OmegaF}\Omega_{t|t}^{i,F}\leq \overline \omega I} where have exploited the definition in (\ref{eq:OmegaFusedwithTilde}) and the facts that $\Omega_{t|t}^i\leq \overline \omega I$ and $\tilde \Omega_t^i\leq  \overline \omega I$, by the transmission rule (\ref{eq:eventTriggered}). We now prove the uniform lower bound. First, notice that if $\Psi_{0|-1}^i>0$ $\forall \, i\in\Nc$ then $\Omega_{t|t}^i>0$ $\forall \, i\in\Nc$. Indeed, by induction we have that if $\Psi_{t|t-1}^i>0$, by the correction step in Algorithm \ref{alg:DKFeventTriggered} it follows that $\Omega_{t|t}^i\geq \Psi_{t|t-1}^i>0$ $\forall \, i\in\Nc$; thus, $\Omega_{t|t}^{i,F}>0$ because it is a convex combination of positive definite matrices; then, it follows $\Omega_{t+1|t}^i=(A(\Omega^{i,F}_{t|t})^{-1}A\tp+Q)^{-1}>0$ by (\ref{eq:infoFormPred}); finally, the fact that $\Psi_{t+1|t}^i>0$ is guaranteed by the choice of $\theta_{t}^i$ as the solution of $\gamma(\Omega_{t+1|t}^i,\theta_t^i)=b$, see \cite{ROBUST_STATE_SPACE_LEVY_NIKOUKHAH_2013,ZORZI_CONTRACTION_CDC}. Let 
$\ones(i)$ denote the indicator function taking value 1 if $i\in\Sc$ and 0 otherwise.
By Lemma \ref{lemma_su_gamma} in Appendix, there exists a constant $\mu>0$ such that 
\al{\Omega_{t|t}^i&=\Psi_{t|t-1}^i +\ones (i)(C^i)\tp(R^i)^{-1}C^i\nn\\
&=\Omega_{t|t-1}^i -\theta_{t-1}^iI+\ones (i)(C^i)\tp(R^i)^{-1}C^i\nn\\
&\geq \mu\Omega_{t|t-1}^i +\ones (i)(C^i)\tp(R^i)^{-1}C^i\nn.}
Moreover, in view of (\ref{upper_OmegaF}), by \cite[Lemma 1 -fact (ii)]{BATTISTELLI2014707} there exists a constant $\upsilon>0$ such that
\al{\label{prima_ineq}\Omega_{t|t}^i&\geq \mu\upsilon A^{-\top}\Omega_{t-1|t-1}^{i,F}A^{-1} +\ones (i)(C^i)\tp(R^i)^{-1}C^i.}
In view of (\ref{eq:eventTriggered}) and (\ref{eq:OmegaFusedwithTilde}), it is not difficult to see that
\al{\label{altra_ineq}\Omega_{t|t}^{i,F}\geq  \frac{1}{(1+\beta)(1+\delta)} \sum_{j\in\Nc} \pi_{i,j} \Omega^j_{t|t}.} Taking into account (\ref{prima_ineq}), we obtain
\al{\Omega_{t|t}^i& \geq \nu A^{-\top}\sum_{j\in\Nc} \pi_{i,j} \Omega^j_{t|t} A^{-1} +\phi_i\nn}
where $\nu:=\mu\upsilon ((1+\beta)(1+\delta))^{-1}>0$ and $\phi_i:=\ones (i)(C^i)\tp(R^i)^{-1}C^i$. The last inequality is similar to the one obtained in the proof of Lemma 1 in \cite{BATTISTELLI201875}. Accordingly, using the same reasonings it is possible to conclude that there exists a constant $\underline \omega>0$ such that $\Omega_{t|t}^i\geq \underline \omega I$. It is also worth noting we have, by (\ref{altra_ineq}),  that $\Omega_{t|t}^{i,F}\geq \underline\omega_F I$ with $\underline\omega_F:=\underline\omega((1+\beta)(1+\delta))^{-1}$.
\end{proof} \smallskip

\prop \label{prop_inequality_e_t}Assume that \textbf{A2}-\textbf{A4} hold. Let  $(q^i_{t|t},\Omega^i_{t|t})$ be the sequence generated according to Algorithm \ref{alg:DKFeventTriggered} with $\Psi_{0|-1}^i>0$  $\forall\, i\in\Nc$. Then, under the least favorable model in (\ref{eq:nominalStateEqgloblf})-(\ref{eq:nominalMeasEqgloblf}),  we have 
\al{\|  e_{t+1}^i\|^2_{\Omega_{t+1|t+1}^i}&\leq\gamma^2 \left( \pi_{i,i}\|e^i_t+\xi_t^i\|^2_{\Omega_{t|t}^i} \right. \nn\\ &\left.+\sum_{j\in\Nc_i} \pi_{i,j}\|e_t^j +\xi_t^i+\eta_t^j\|^2_{\Omega_{t|t}^j}\right),\quad \forall \,  i\in\Nc\nn}
where $0<\gamma<1$ is a constant and 
\al{ 
    \xi_{t}^{i} &:= A^{-1}[B\tilde w_{t}-\mathbf{1}_{\mathcal{S}}(i)(\Psi_{t+1|t}^{i})^{-1}(C^{i})\tp(R^{i})^{-1}D^i\tilde v_{t+1}^{i}],\nn\\
    \eta_{t}^{j} &:= (1-c^{j}_{t}) (\hat{x}^{j}_{t|t}-\bar{x}^{j}_{t}),\quad \bar x_{t}^j:=(\bar \Psi_t^j)^{-1} \bar q_t^j.	\nn}
\eprop

\begin{proof} Recall that $q_{t+1|t}^i=\Psi^i_{t+1|t}x_{t+1|t}^i$ and $q_{t|t}^i=\Omega^i_{t|t}x_{t|t}^i$. Then, it is not difficult to see that
 \begin{align*}
   & {x}_{t+1|t+1}^{i}= (\Omega_{t+1|t+1}^{i})^{-1}  \left ( \Psi_{t+1|t}^{i}{x}_{t+1|t}^{i}\right.\\&\hspace{0.2cm}\left.+ \mathbf{1}_{\mathcal{S}}(i)(C^{i})\tp(R^{i})^{-1}C^{i}x_{t+1}+\mathbf{1}_{\mathcal{S}}(i)(C^{i})\tp(R^{i})^{-1}D^i\tilde v_{t+1}^{i}\right )\\
         &  x_{t+1}  = (\Omega_{t+1|t+1}^{i})^{-1}\left ( \Psi_{t+1|t}^{i}x_{t+1}\right.\\ &\hspace{0.2cm}+\left.\mathbf{1}_{\mathcal{S}}(i)(C^{i})\tp(R^{i})^{-1}C^{i}x_{t+1}\right )
\end{align*}
which implies
\begin{align*}
    e_{t+1}^{i} &= x_{t+1}-{x}_{t+1|t+1}^{i}\nn\\ &= (\Omega_{t+1|t+1}^{i})^{-1}\Psi_{t+1|t}^{i}(x_{t+1}- {x}_{t+1|t}^{i}+\bar{v}_{t+1}^{i})
\end{align*}
where $\bar{v}_{t+1}^{i} = -\mathbf{1}_{\mathcal{S}}(i)(\Psi_{t+1|t}^{i})^{-1}(C^{i})\tp(R^{i})^{-1}D^i\tilde v_{t+1}$. Since $\Psi^i_{t+1|t}\leq \Omega^i_{t+1|t+1}$, we have 
\begin{align}
    \|&e_{t+1}^{i}\|_{\Omega_{t+1|t+1}^{i}}^{2} \nn\\ &= \|x_{t+1}-{x}_{t+1|t}^{i}+\bar{v}_{t+1}^i \|^2_{\Psi^i_{t+1|t} (\Omega_{t+1|t+1}^{i})^{-1}\Psi^i_{t+1|t}} \nn\\
          &\leq \|x_{t+1}-{x}_{t+1|t}^{i}+\bar{v}_{t+1}^i\|_{\Psi_{t+1|t}^{i}}^{2}. \label{eq:lemma2_22}
\end{align}
By Lemma 1 - fact (iii) in \cite{BATTISTELLI2014707} there exists a constant $0<\gamma\leq 1$ such that \begin{equation*}
   \Psi_{t+1|t}^i= \Omega_{t+1|t}^{i}-\theta_t^i \leq \Omega_{t+1|t}^{i}\leq \gamma^{2}A^{-T}\Omega_{t|t}^{i,F}A^{-1}
\end{equation*}
where we exploited the fact that $\Omega^{i,F}_{t|t} \geq \underline \omega_F I$ by Proposition \ref{prop_bound_Omega}. 
Taking into account (\ref{eq:lemma2_22}), it follows that
\begin{align}
    \|e_{t+1}^{i}\|_{\Omega_{t+1|t+1}^{i}}^{2} &\leq \gamma^2\|x_{t+1}-{x}_{t+1|t}^{i}+\bar{v}_{t+1}^i\|_{A^{-\top}\Omega_{t|t}^{i,F}A^{-1}}^{2}\nn \\
    &=\gamma^2 \|A(x_{t}-{x}^{i,F}_{t|t})+B \tilde w_{t}+ \bar{v}^{i}_{t+1}\|_{A^{-\top}\Omega_{t|t}^{i,F}A^{-1}}^{2} \nn\\
 &= \gamma^{2}\|x_{t}- {x}^{i,F}_{t|t} + \xi^{i}_{t}\|_{\Omega_{t|t}^{i,F}}^{2}  \nn
\end{align}
where we exploited the fact that ${x}^{i,F}_{t|t} = (\Omega_{t|t}^{i,F})^{-1}q^{i,F}_{t|t}$ and $x_{t+1|t}^i=Ax_{t|t}^{i,F}$. Then, using the same reasonings in the second part of the proof of Lemma 2 in \cite{BATTISTELLI201875} it is not difficult to see that 
{\small \begin{align*} 
    \|e_{t+1}^{i}  \|_{\Omega_{t+1|t+1}^{i}}^{2}& \leq \gamma^{2}\left ( \pi_{i,i}\|{e}^{i}_{t}+ \xi^{i}_{t}\|_{\Omega_{t|t}^{i}}^{2} \right.\\
    &\left.+\sum_{j\in\mathcal{N}_{i}}\pi_{i,j}\|{e}^{j}_{t} + \xi^{i}_{t}+\eta_t^j\|_{\Omega_{t|t}^{j}}^{2} \right )
\end{align*}} where
\al{&\eta_t^j:= x_{t|t}^j-\check x_{t|t}^j,\quad \check x_{t|t}^j=(\check \Omega_{t|t}^j)^{-1}\check q_{t|t}^j\nn\\ & \check{q}^{j}_{t|t}:= c_{t}^{j}q^{j}_{t|t}+(1-c_{t}^{j})\tilde{q}^{j}_{t}, \quad 
\check{\Omega}^{j}_{t|t}:= c_{t}^{j}\Omega^{j}_{t|t}+(1-c_{t}^{j})\tilde{\Omega}^{j}_{t}\nn}
which concludes the proof.
\end{proof}
\smallskip
In view of Proposition \ref{propo_leastfav} and Proposition \ref{prop_bound_Omega}, it follows that $\xi_t^i$, with $i\in\Nc$, is uniformly bounded in mean-square under the assumptions \textbf{A1}-\textbf{A4}. Indeed, if $i\notin\Sc$, then
\al{\xi_t^i&=A^{-1} B \tilde w_t=A^{-1} \left[\begin{array}{cc}
 B & 0\end{array}\right]\left[ \begin{array}{c}
\tilde w_t \\ \tilde v_{t+1}\end{array}\right]\nn\\
\Var(\xi_t^i)&=A^{-1}\left[\begin{array}{cc}
 B & 0\end{array}\right]\Var\left(\left[ \begin{array}{c}
\tilde w_t \\ \tilde v_{t+1}\end{array}\right]\right)\left[\begin{array}{c}
 B\tp \\ 0\end{array}\right] A^{-\top}\nn\\
&\leq \overline \rho  \tp A^{-1}BB\tp A^{-\top}\leq \overline \rho \sigma_{max}( A^{-1}BB\tp A^{-\top})I;\nn}
if $i\in\Sc$, then
\al{\xi_t^i&=A^{-1}\left[\begin{array}{cc}
 B& H_i\end{array}\right]\left[ \begin{array}{c}
\tilde w_t \\ \tilde v_{t+1}\end{array}\right]\nn} where 
\al{H_i:=[\,
 0\; \ldots \;\underbrace{-(\Psi^i_{t+1|t})^{-1}(C^i)\tp(R^i)^{-1}D^i}_{\hbox{$i+1$-th block}} \;\ldots \; 0\,]\nn}
 and thus
\al{\Var(\xi_t^i)&=A^{-1}\left[\begin{array}{cc}
 B& H_i\end{array}\right]\Var\left(\left[ \begin{array}{c}
\tilde w_t \\ \tilde v_{t+1}\end{array}\right]\right)\left[ \begin{array}{c}
B\tp \\ H_i\tp\end{array}\right] A^{-\top}\nn\\
&\leq \overline \rho A^{-1}(BB\tp+H_iH_i\tp )A^{-\top}\nn\\ &\leq  \overline \rho\sigma_{max}(A^{-1}(BB\tp+H_iH_i\tp )A^{-\top})I.\nn}
We conclude that $\Var[\xi_t^i]\leq \rho_\xi^2$ where 
\al{\rho_\xi^2:=\overline \rho\,\underset{i\in\Nc}{\max}\,\sigma_{max}(A^{-1}(BB\tp+H_iH_i\tp )A^{-\top}).\nn}
Finally, it is worth noting that $\|\eta_t^j\|^2_{\Omega_{t|t}^j}\leq \alpha$, indeed if $c_t^j=0$, i.e. there is not transmission from node $j$, then $\|\eta_t^j\|^2_{\Omega_{t|t}^j}=\|x_{t|t}^j-\bar x_t^j\|^2_{\Omega_{t|t}^j}\leq \alpha$ by the transmission rule in (\ref{eq:eventTriggered}); otherwise, if $c_t^i=1$ then $\|\eta_t^j\|^2_{\Omega_{t|t}^j}=0\leq \alpha$.
\smallskip
\teo \label{teo_conv}Assume that the hypotheses \textbf{A1}-\textbf{A4} hold. Then, the estimation error $e_t^i$, with $i\in\Nc$, is uniformly bounded in mean-square as $t\rightarrow \infty$, i.e.
\al{\lim_{t\rightarrow \infty} \Es[\|e_t^i\|^2]\leq \left(\frac{\sqrt 2\gamma}{1-\gamma} \frac{\sum_{i\in\Nc} \sqrt{p_i}(\sqrt{\overline \omega} \rho_\xi +\sqrt{\alpha}) }{\sqrt{\underline \omega}\min_{i\in\Nc} \sqrt{p_i}}\right)^2\nn}
where $\underline \omega, \overline\omega,\alpha,\gamma$ are the constants given in Proposition \ref{prop_bound_Omega} and Proposition \ref{prop_inequality_e_t}.
\eteo
\begin{proof} Using the same reasonings in the first part of the proof of Theorem 1 in \cite{BATTISTELLI201875}, it is not difficult to prove that
\al{&\label{ineq1_teo}\sqrt{\Es[\Vc_{t+1}(e_{t+1})]}\leq \gamma\left(\sqrt{\Es[\Vc_t(e_t)]}+\sqrt{\Es[\Vc_t(\omega_t)]}\right)\\ 
&\label{ineq2_teo}\sqrt{\Vc_t(\omega_t)}\leq \sum_{i\in\Nc}
\sqrt{p_i}  \sqrt{\Es[\|\omega_t^i\|^2_{\Omega_{t|t}^i} }   }
where $0<\gamma<1$ is the constant given by Proposition \ref{prop_inequality_e_t}, $\omega_t:=\col{\omega_t^i;i\in\Nc}$,
\al{\omega_t^i:=\left\{\begin{array}{cc} \bar \omega_t^i, & \hbox{ if } e_t^i =0\\
 \|\bar \omega_t^i\|^2_{\Omega_{t|t}^i}/ \|e_t^i\|^2_{\Omega_{t|t}^i} e_t^i, & \hbox{otherwise} \end{array}\right.\nn} 
 and $\bar \omega_t^i$ is the vector in the set $\{\xi_t^i\}\cup \{\xi_t^j+\eta_t^i,\, j\in\overline\Nc_i\}$ maximizing the weighted norm $\|\cdot \|_{\Omega_{t|t}^i}$. Notice that
{\small  \al{0\leq \Es[\|\xi_t^j-\eta_t^i\|^2_{\Omega_{t|t}^i}]
 =\Es[\|\xi_t^j \|^2_{\Omega_{t|t}^i} + \| \eta_t^i\|^2_{\Omega_{t|t}^i}-2(\xi_t^j)\tp\Omega_{t|t}^i \eta_t^i]
\nn}}
and thus 
\al{\Es[2(\xi_t^j)\tp\Omega_{t|t}^i \eta_t^i] \leq \Es[\|\xi_t^j \|^2_{\Omega_{t|t}^i} + \| \eta_t^i\|^2_{\Omega_{t|t}^i}].\nn}
Moreover, 
{\small \al{ &\Es  [\|\omega_t^i\|^2_{\Omega_{t|t}^i}]=\Es[\|\bar \omega_t^i\|^2_{\Omega_{t|t}^i}]\nn\\ &\leq 
\max\{\Es[\|\xi_t^i\|^2_{\Omega_{t|t}^i}],\max_{j\in\Nc}\Es[\| \xi_t^j\|^2_{\Omega_{t|t}^i}+\| \eta_t^i\|^2_{\Omega_{t|t}^i}+2(\xi_t^j)\tp\Omega_{t|t}^i \eta_t^i]\}
\nn\\
&\leq 
\max\{\Es[\|\xi_t^i\|^2_{\Omega_{t|t}^i}],\max_{j\in\Nc}2\,\Es[\| \xi_t^j\|^2_{\Omega_{t|t}^i}+\| \eta_t^i\|^2_{\Omega_{t|t}^i}]\}\nn \\
&\leq 2\,
\max_{j\in\Nc}\Es[\| \xi_t^j\|^2_{\Omega_{t|t}^i}+\| \eta_t^i\|^2_{\Omega_{t|t}^i}].\nn 
} } Taking into account (\ref{ineq2_teo}), we obtain
\al{\sqrt{\Vc_t(\omega_t)}&\leq \sqrt{2}\sum_{i\in\Nc}
\sqrt{p_i}  \max_{j\in\Nc} \sqrt{ \Es[\| \xi_t^j\|^2_{\Omega_{t|t}^i}+\| \eta_t^i\|^2_{\Omega_{t|t}^i}] }\nn\\
&\leq \sqrt{2}\sum_{i\in\Nc}
\sqrt{p_i}  \left(\max_{j\in\Nc} \sqrt{ \Es[\| \xi_t^j\|^2_{\Omega_{t|t}^i}}+\sqrt{\Es[\| \eta_t^i\|^2_{\Omega_{t|t}^i}] }\right)\nn\\ 
&\leq \sqrt{2}\sum_{i\in\Nc}
\sqrt{p_i}  \left(\sqrt{\overline\omega}\rho_\xi+ \sqrt{\alpha}\right)\nn}
where we exploited the upper bounds for $\|\xi_t^i\|^2_{\Omega_{t|t}^j}$  and $\|\eta_t^i\|^2_{\Omega_{t|t}^i}$. Taking into account (\ref{ineq1_teo}) we obtain 
\al{&\sqrt{\Es[\Vc_{t+1}(e_{t+1})]}\nn\\ &\hspace{0.8cm}\leq \gamma  \sqrt{\Es[\Vc_t(e_t)]}+ \sqrt{2}\gamma\sum_{i\in\Nc}
\sqrt{p_i}  \left(\sqrt{\overline\omega}\rho_\xi+ \sqrt{\alpha}\right) \nn}
which implies 
\al{\underset{t\rightarrow \infty}{\mathrm{lim\,\sup}} \,\sqrt{\Es[\Vc_{t}(e_{t})]}\leq  \frac{\sqrt 2\gamma}{1-\gamma} \sum_{i\in\Nc}
\sqrt{p_i}  \left(\sqrt{\overline\omega}\rho_\xi+ \sqrt{\alpha}\right)\nn.} Finally, recalling that $\Omega_{t|t}^i\geq \underline \omega I$ by Proposition \ref{prop_bound_Omega},  it is sufficient to note that 
\al{\Es[\|e_t\|^2]\leq \frac{\Es[\Vc_t(e_t)]}{ \underline \omega \min_{i\in\Nc} p_i}\nn}
to obtain the claim. 
\end{proof}


 \section{A robust strategy with local tolerances}\label{sec:local_algo}
In Section \ref{sec:algo} we have proposed a distributed strategy with event-triggered communication which provides a solution to the minimax problem in (\ref{central_pb}) which is suboptimal. Indeed, at each sensor node the state prediction in (\ref{eq:infoFormPred}) is the solution to the (local) minimax problem 
\al{ x_{t+1|t}^i=\underset{g_t \in\mathcal G_i}{\argmin}\underset{\tilde \phi_t^i\in\mathcal B_t^i}{\max}\, \tilde{\mathbb E}_i[\| x_{t+1}-g_t(y_t^i)\|^2|Y_{t-1}]\label{local_pb}}
where 
\begin{align}
    \tilde{\mathbb{E}}_i&\left [ \|x_{t+1}-g_{t}(y_{t}^i)\|^{2}|Y_{t-1} \right]: =\nn\\ & \int_{\mathbb R^{n}}\int_{\mathbb R^{n+p_i}}\|x_{t+1}-g_{t}(y_{t})\|^{2}\tilde{\phi}_{t}^i(z_{t}^i|x_{t})\check{f}_{t}(x_{t}|Y_{t-1})dz_{t}^idx_{t}\nn
\end{align} with $z_t^i:=[\, x_{t+1}^\top \; (y_t^i)^\top\,]^\top$;  $\Bc_t^i$ is the ambiguity set about the local and nominal model in (\ref{eq:nominalStateEq})-(\ref{eq:nominalMeasEq}), whose nominal density is denoted by $\phi_t^i$, with  tolerance $b$; $\mathcal G_i$ is the set of estimators having finite second order moments for any $\tilde \phi_t^i\in \mathcal B_t^i$. Thus, the local least favorable transition probability density $\tilde \phi_t^i$ solution to (\ref{local_pb}) does not necessarily agree with the global least favorable density $\tilde \phi_t$ obtained solving the centralized problem in (\ref{central_pb}). 

Let 
\al{\bar f_t(z_t|Y_{t-1}):=&\int_{\Rs^{n}} \phi_t(z_t|x_t) \check f_t(x_t|Y_{t-1})  \mathrm dx_t\nn\\
\tilde f_t(z_t|Y_{t-1}):=&\int_{\Rs^{n}}\tilde \phi_t(z_t|x_t) \check f_t(x_t|Y_{t-1})  \mathrm dx_t\nn
}
be the pseudo-nominal and the least favorable conditional probability densities of $z_t$ given $Y_{t-1}$, respectively. In a similar way we define the  pseudo-nominal and the least favorable conditional probability densities of $z_t^i$ given $Y_{t-1}$, respectively, as

\al{\bar f_t^i(z_t|Y_{t-1}):=&\int_{\Rs^{n }} \phi_t^i(z_t|x_t) \check f_t(x_t|Y_{t-1})  \mathrm dx_t\nn\\
\tilde f_t^i(z_t|Y_{t-1}):=&\int_{\Rs^{n }} \tilde \phi_t^i(z_t|x_t) \check f_t(x_t|Y_{t-1})  \mathrm dx_t.\nn
} In \cite{ROBUST_STATE_SPACE_LEVY_NIKOUKHAH_2013} it has been shown that the least favorable density $\tilde \phi^i_t$ solution to (\ref{local_pb}) is such that the Kullback-Leibler divergence between $\bar f_t^i$ and $\tilde f_t^i$, i.e. 
\al{D_{KL}(\tilde f_t^i,\bar f_t^i):=\int_{\Rs^{n+p_i}}\tilde f_t^{i}(z_t^i|Y_{t-1}) \log\left(\frac{\tilde f_t^{i}(z_t^i|Y_{t-1})}{\bar f_t^{i}(z_t^i|Y_{t-1})}\right)\mathrm d z_t^i\nn}
is equal to the tolerance of the ambiguity set $\Bc_t^i$ used in (\ref{local_pb}).

Drawing inspiration from \cite{s20113244}, it is possible to mitigate the fact that the local least favorable density $\tilde \phi_t^{i}$ is different from the one obtained from the global least favorable density $\tilde \phi_t$ by considering the ambiguity set $\Bc_t^i$ about the nominal density $\phi^i_t$ and with ``local'' tolerance \al{\label{b_time_var}b_t^i&=D_{KL}(\tilde f^{i,g}_t, \bar f^i_t):=\nn\\ & \int_{\mathbb R^{n+p_i}} \tilde f_t^{i,g}(z_t^i|Y_{t-1})\log\left( \frac{\tilde f_t^{i,g}(z_t^i|Y_{t-1})}{\bar f_t^{i }(z_t^i|Y_{t-1})}\right)\mathrm d z_t^i }  where $\tilde f_t^{i,g}$ is the least favorable density at the sensor node $i$ obtained by marginalizing $\tilde f_t$ with respect to $y_t^l$ with $l\neq i$. Although such a choice does not guarantee that $\tilde \phi_t^i$, i.e. the least favorable density solution to (\ref{local_pb}) with tolerance (\ref{b_time_var}), leads to a density $\tilde f_t^i$ which coincides with $\tilde f_t^{i,g}$, at least it is ensured that
\al{ D_{KL}(\tilde f^{i}_t,\bar f_t^i)=D_{KL}(\tilde f_t^{i,g},\bar f_t^i).\nn}

In \cite{zorzi2018lf} it has been shown that if $b>0$ in (\ref{central_pb}) is taken sufficiently small, then the least favorable model $\tilde \phi_t$ solution to (\ref{central_pb}) converges to a stationary Gaussian model as $t\rightarrow \infty$. Moreover, the conditional mean of $z_t$ given $Y_{t-1}$ under $\tilde f_t$ is the same of the one under $\bar f_t$, see \cite[Theorem 1]{s20113244}. Let $K$ and $\tilde K$ denote the asymptotic covariance matrices of $z_t$ given $Y_{t-1}$ under the Gaussian densities $\bar f_t$ and $\tilde f_t$, respectively. Then, $b_t^i\rightarrow b^i$ as $t\rightarrow \infty$ where 
\al{ b^i=
      \frac{1}{2}\left[ \log \det (K_{{i}}\tilde{K}_{{i}}^{-1})+\tr(\tilde{K}_{{i}}K^{-1}_{{i}}) -(n+p_{i}) \right];\nn}
$K_i$ is the matrix obtained from $K$ by deleting the rows and the columns corresponding to $y_t^l$ with $l\neq i$; $\tilde K_i$ is obtained from $\tilde K$ likewise.

In conclusion, we propose the following distributed filtering strategy with event-triggered communication. The central unit can compute offline the asymptotic tolerance $b^i$ for each sensor node $i\in \Sc$ from the global least favorable density $\tilde \phi_t$ and send it to the corresponding sensor node. After this offline step, the remaining part of the algorithm is as Algorithm \ref{alg:DKFeventTriggered} with the exception of the second equation in (\ref{eq:infoFormPred}) and (\ref{eq:estPredictionStep}) which become, respectively,
\al{\text{Find } &\theta_{t}^i>0 \text{ s.t. }\gamma(\Omega_{t+1|t}^i, \theta_{t}^i)=b^i \nn\\
                \text{Find } &\bar\theta_{t}^i>0 \text{ s.t. }\gamma(\bar \Omega_{t+1}^i, \bar\theta_{t}^i)=b^i.\nn}

In regard to the stability properties of this algorithm with local tolerances, it is not difficult to see that under Assumptions \textbf{A1}-\textbf{A4}
\al{\lim_{t\rightarrow \infty} \Es[\|e_t^i\|^2]\leq \left(\frac{\sqrt 2\gamma}{1-\gamma} \frac{\sum_{i\in\Nc} \sqrt{p_i}(\sqrt{\overline \omega} \rho_\xi +\sqrt{\alpha}) }{\sqrt{\underline \omega}\min_{i\in\Nc} \sqrt{p_i}}\right)^2\nn} for any $i\in\Nc$. 
The proof follows the same ideas exploited in Section \ref{sec_stab}. The unique difference is the derivation of the constant $\underline \omega$. More precisely, in the proof of Proposition \ref{prop_bound_Omega} the constant $\mu$ is derived as follows. Since each sensor node $i$ has its own tolerance $b^i$, then we have that 
$\gamma(\Omega^i_{t+1|t},\theta_t^i)=b^i$. By Lemma \ref{lemma_su_gamma} in Appendix, it follows that there exists $\mu_i>0$ such that $\Omega_{t|t}^i\geq \mu_i\Omega_{t|t-1}^i$ $\forall \,i\in\Sc$. Then, $\Omega_{t|t}^i\geq \mu \Omega_{t|t-1}^i$ $\forall \,i\in\Sc$ with $\mu:=\min_{i\in\Sc} \mu_i>0$.  
 \section{Simulations}\label{sec:numerical}
 
In this section we evaluate the performance of the proposed robust distributed Kalman algorithms with event-triggered communication. We consider the problem of tracking the position of a target by using noisy position measurements obtained by a network of $N=100$ nodes where $20$ of them are sensor nodes. The possible connections among the nodes has been randomly generated in such a way the network is strongly connected. Moreover,  the percentage of node connections is equal to $4\%$. The model for the motion of the target is 
\al{\label{proj_cont}\dot x^c_s =\Phi x^c_s+\dot w^c_s,\quad s\in\Rs}
where 
\al{\Phi=\left[\begin{array}{cc}0 & 0 \\ \mathbf I_3 & 0\end{array}\right],\nn}
$x^c_s=[\, v_{x,s}\;v_{y,s}\;v_{z,s}\;p_{x,s}\;p_{y,s}\;p_{z,s}\;\,]\tp$ with $v$ denoting the velocity, $p$ the position and the subscripts $x,y,z$ denoting the three spatial dimensions; $w_s^c$ is a Wiener process with zero mean and rate of variance equal to 0.1. We discretize (\ref{proj_cont}) with sampling time equal to $0.1$. The corresponding discrete time model is $x_{t+1}=Ax_{t}+Bw_t$ where $x_t$ is the sampled version of $x^c_s$, $A=\mathbf I_6+0.1\Phi$, $B=\sqrt{0.001}I$, $w_t$ is a zero-mean normalized Gaussian white noise and thus $Bw_t$ is the sampled version of $w_s^c$. We assume that every sensor measures the position of the target in either two horizontal dimensions, or a combination of one horizontal dimension and the vertical dimension; in plain words, one sensor does not have measurements in all the three dimensions. Therefore, we obtain the  nominal discrete state-space model (\ref{eq:nominalStateEq})-(\ref{eq:nominalMeasEq}) where $C^i=[\, 0\; 0 \;0 \;\mathrm{diag(1,1,0)}\,]$, in the case that the sensor measures only the horizontal positions, or $C^i=[\, 0\; 0 \;0 \;\mathrm{diag(1,0,1)}\,]$, $C^i=[\, 0\; 0 \;0 \;\mathrm{diag(0,1,1)}\,]$, in the case that the sensor measures one horizontal position and the vertical position. Moreover,  $R^i=D^i(D^i)\tp=\sqrt{k} PR_0P\tp$ where $R_0=0.5\cdot \mathrm{diag}(1,4,7)$ and $P$ is a permutation matrix randomly chosen for every node. Finally, the initial state $x_0$ is a Gaussian random vector with zero-mean and covariance matrix $V_{0|-1}=I$. Since the previous model is just an idealization of the underlying physical system, we assume that the actual state-space model belongs to the ambiguity set $\mathcal B_t$ about the aforementioned nominal model and with tolerance 
$b=0.05$. More precisely, we assume that the actual model is  the least favorable model solution to (i.e. the maximizer of) the centralized problem in (\ref{central_pb}).

In the following we consider the distributed algorithms:
\begin{itemize}
\item \textsf{RDKF} -- the distributed robust Kalman filter with event-triggered communication in Algorithm \ref{alg:DKFeventTriggered} with $b=0.05$. Here, the 
transmission rule (\ref{eq:eventTriggered}) is with $\alpha=10$, $\beta=0.2$ and $\delta=0.5$;
\item \textsf{RDKFLOC} --  the distributed robust Kalman filter with event-triggered communication of Section \ref{sec:local_algo} where each sensor nodes $i $ has its own tolerance $b^i$ computed from the global least favorable model; the 
transmission rule (\ref{eq:eventTriggered}) is with $\alpha=10$, $\beta=0.2$ and $\delta=0.5$;
\item \textsf{DKF1} -- the distributed Kalman filter with event-triggered communication proposed in \cite{BATTISTELLI201875} and the transmission rule is with $\alpha=10$, $\beta=0.2$ and $\delta=0.5$;
\item \textsf{DKF2} -- is the same as \textsf{DKF1} but the transmission rule is with $\alpha=0.01$, $\beta=0.2$ and $\delta=0.5$.
\end{itemize}
It is worth noting that \textsf{RDKF}, \textsf{RDKFLOC} and \textsf{DKF1} have the same parameters for the transmission rule. As we will see later, taking the parameter $\alpha$ in (\ref{eq:eventTriggered}) the same for \textsf{RDKF}, \textsf{RDKFLOC} and \textsf{DKF1} provides a transmission rate for \textsf{DKF1} which is smaller than the robust versions. For this reason, we also consider \textsf{DKF2} where the parameter $\alpha$ has been decreased in oder to increase the transmission rate. We also tried to increase the transmission rate by keeping fixed $\alpha=10$ and changing $\beta,\delta$; however, we did not notice a significant growth in terms of transmission rate.  \begin{figure}
\begin{center}
\includegraphics[width=0.5\textwidth]{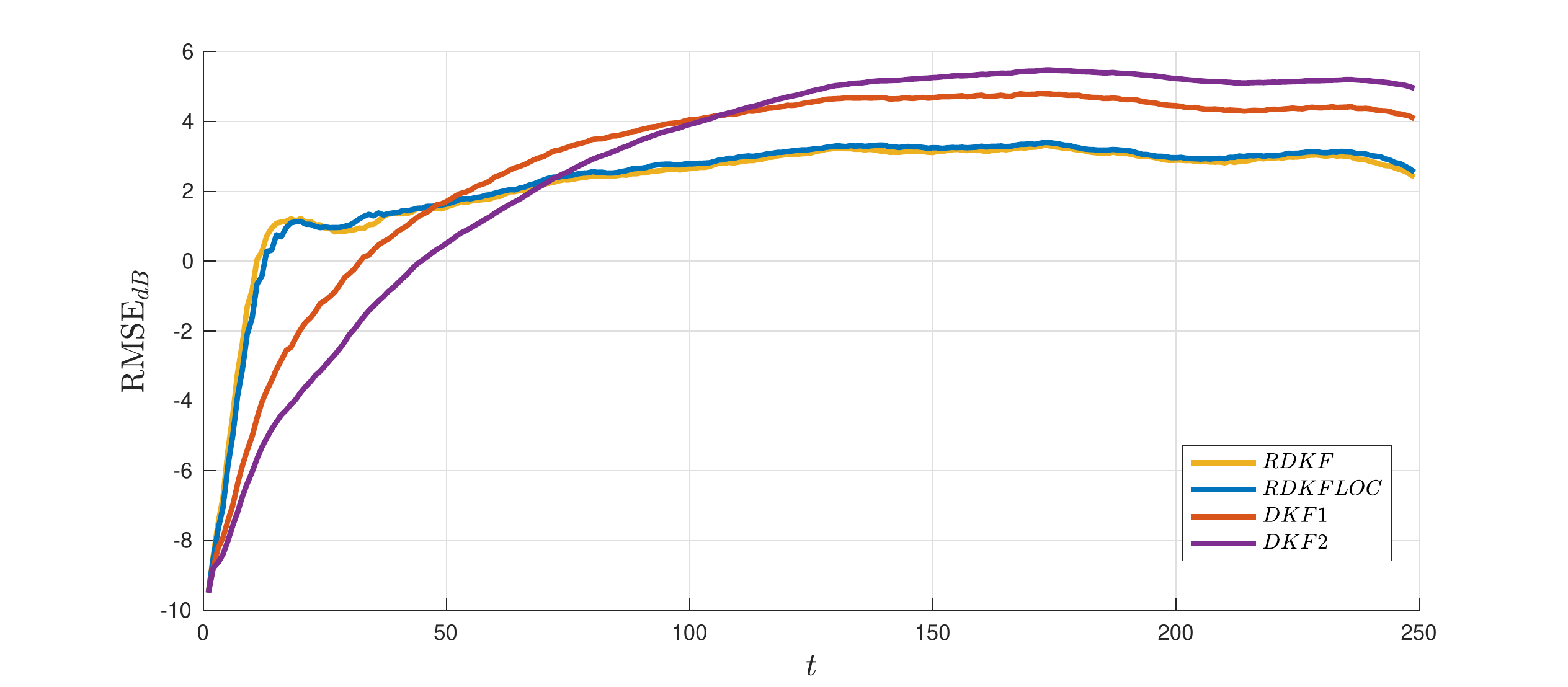}    
\caption{Average RMSE across the network.} 
\label{fig:rmse}
\end{center}
\end{figure}

 \begin{figure*}
\begin{center}
\includegraphics[width=\textwidth]{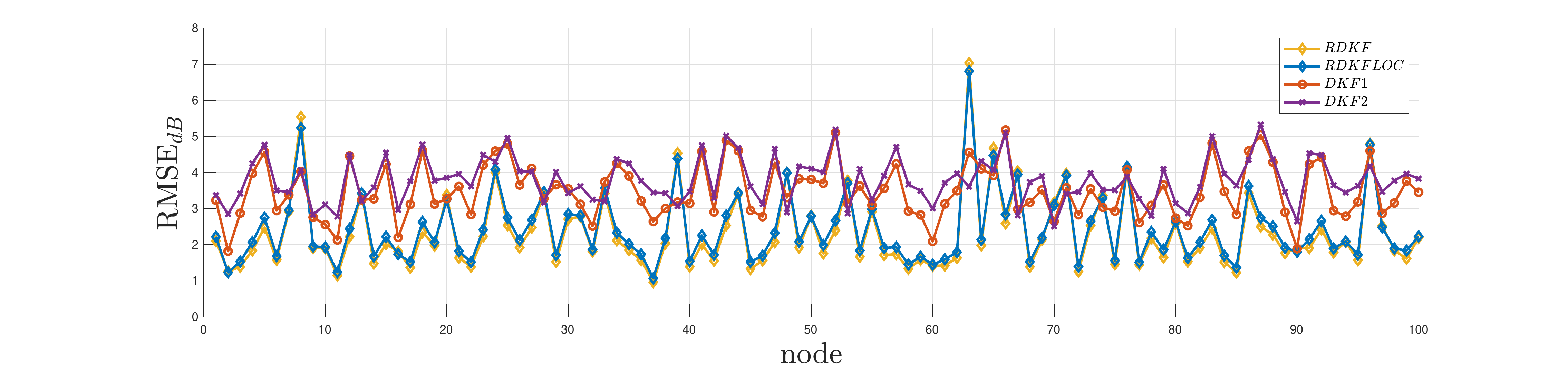}    
\caption{Average RMSE at each node over the time horizon $[1,250]$.} 
\label{fig:rmsen}
\end{center}
\end{figure*}

 \begin{figure}
    \begin{center}
    \includegraphics[width=0.5\textwidth]{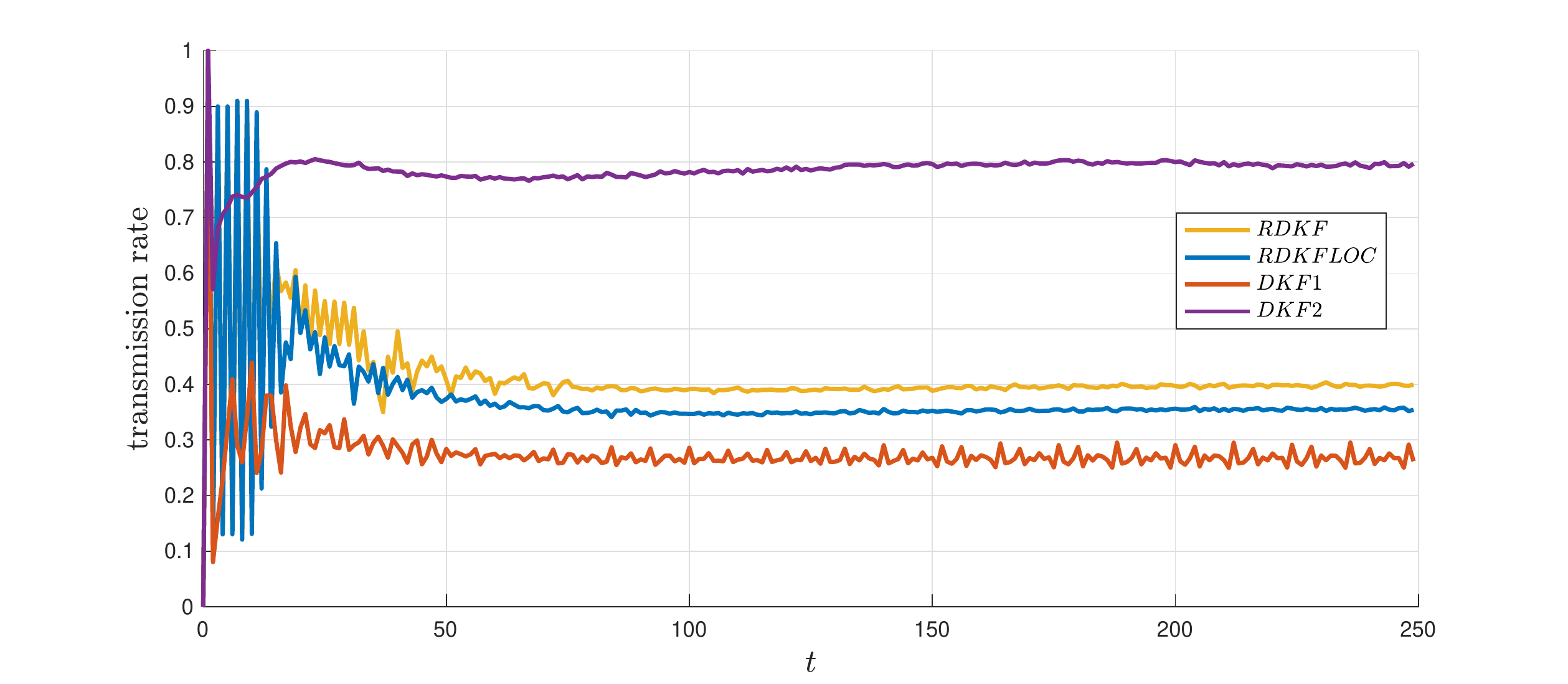}    
    \caption{Average transmission rate across the network.} 
    \label{fig:erre}
    \end{center}
\end{figure}
 \begin{figure}
    \begin{center}
    \includegraphics[width=0.5\textwidth]{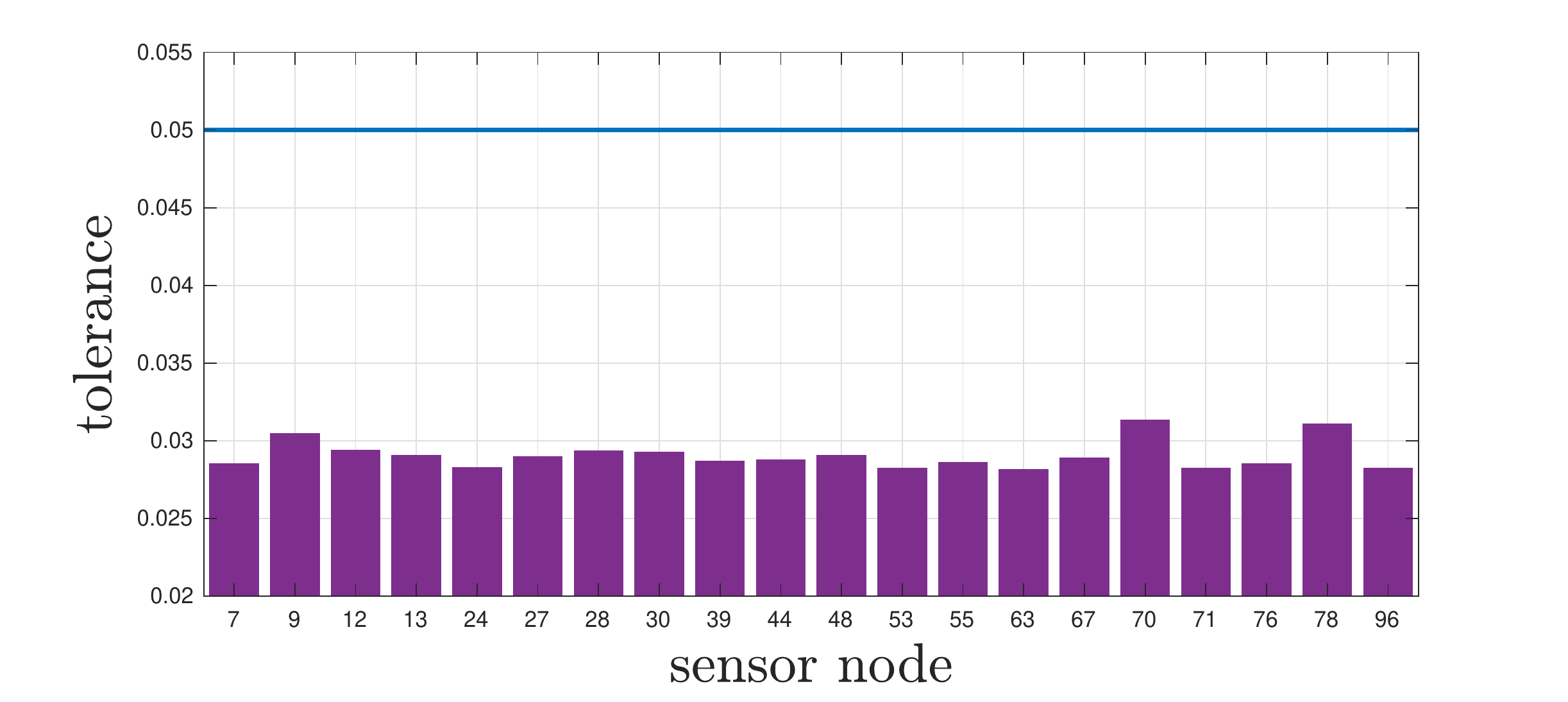}    
    \caption{Tolerance used in \textsf{RDKF} (blue line) versus the local tolerances at each sensors node used in \textsf{RDKFLOC} (purple bars).} 
    \label{fig:tolloc}
    \end{center}
\end{figure}
 \begin{figure}
    \begin{center}
    \includegraphics[width=0.5\textwidth]{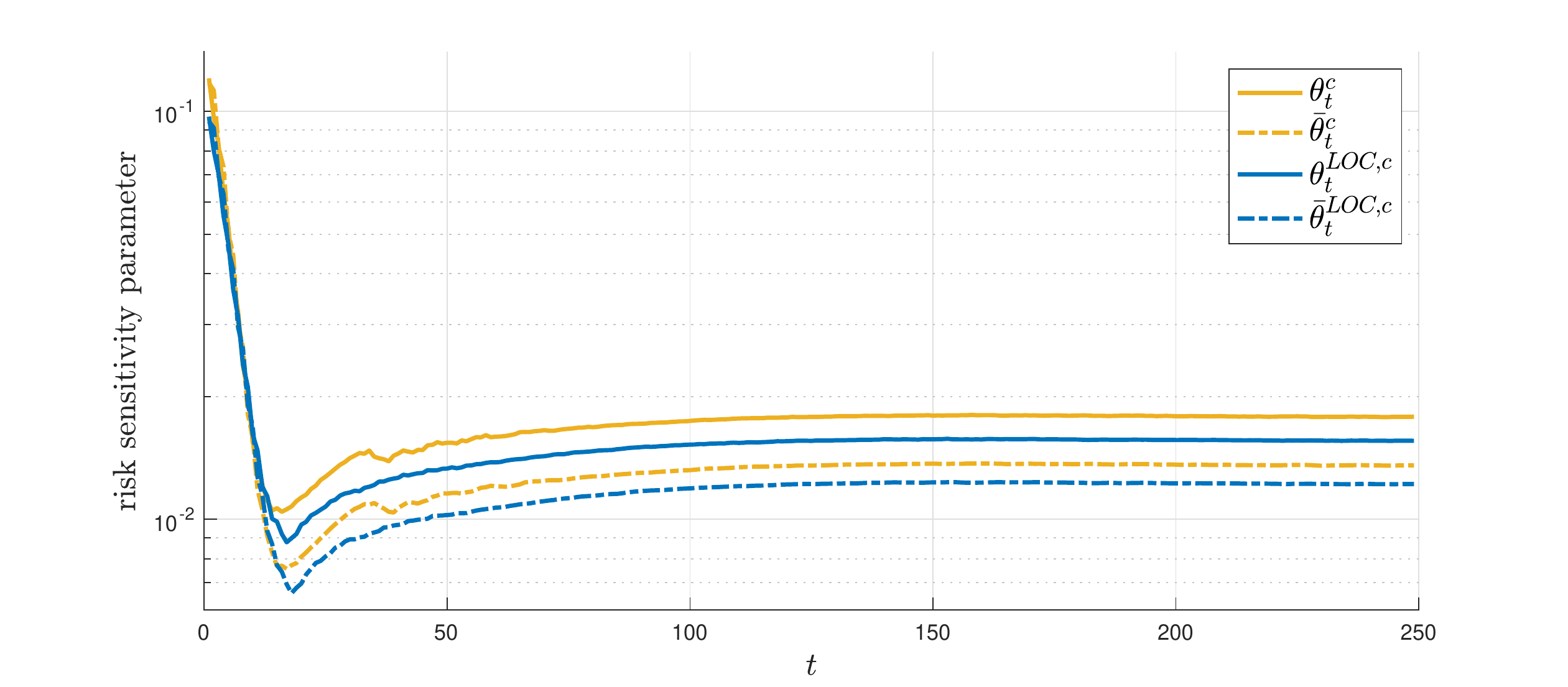}    
    \caption{Average risk sensitivity parameters across  the communication nodes of the network (in logarithmic scale).} 
    \label{fig:thetac}
    \end{center}
\end{figure}
 \begin{figure}
    \begin{center}
    \includegraphics[width=0.5\textwidth]{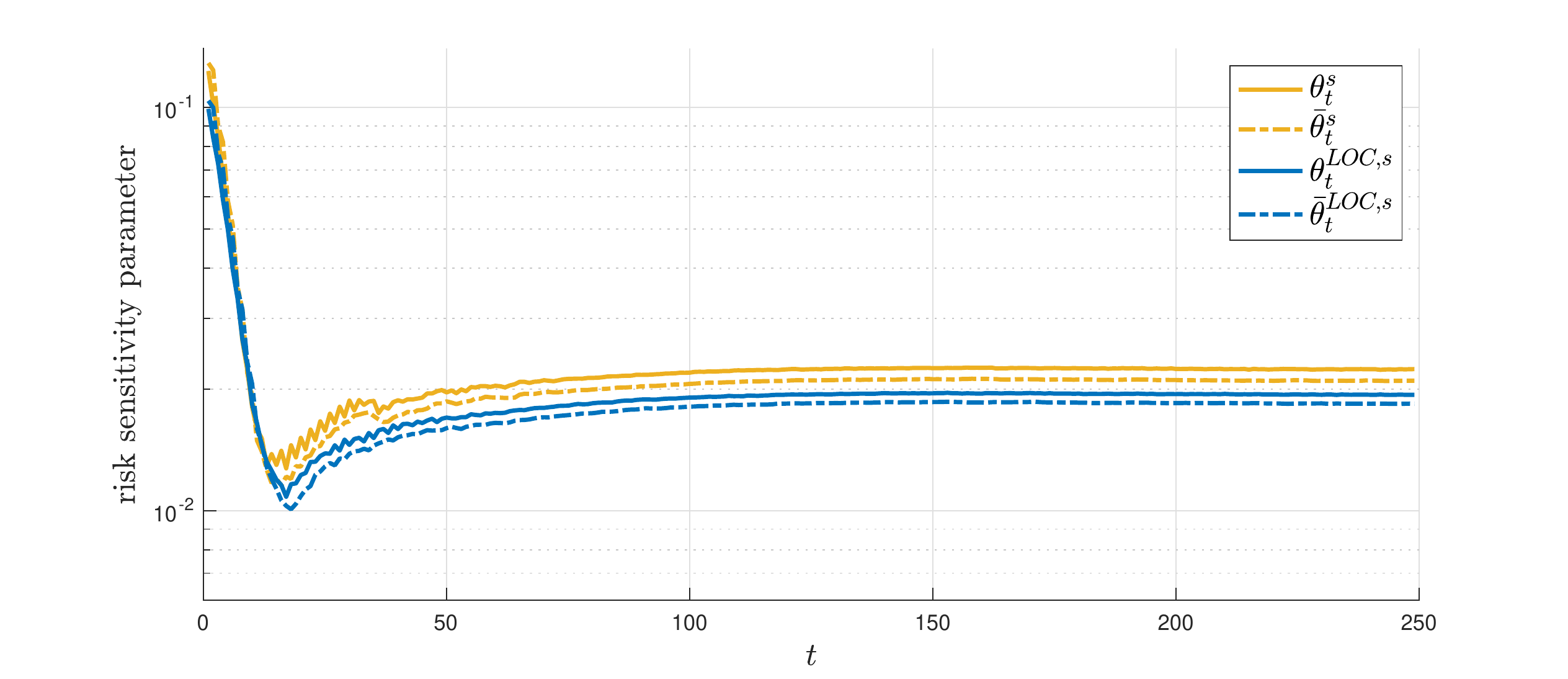}    
    \caption{Average risk sensitivity parameter across the  sensor nodes of the network (in logarithmic scale).} 
    \label{fig:thetas}
    \end{center}
\end{figure}

We consider a Monte Carlo study which is composed by  500 independent runs (which correspond to different target paths generated by the least favorable model \footnote{For more details on how to generate a realization from the least favorable model see \cite[Section V]{ROBUST_STATE_SPACE_LEVY_NIKOUKHAH_2013}.}) over a time horizon of 250 seconds. For each run, we estimate the state using the aforementioned distributed algorithms and for each of them we consider the following performance indexes:
\begin{itemize}
\item The average root mean square error across the network at time $t$: $$ \mathrm{RMSE}_t=\frac{1}{N}\sum_{i=1}^N \|x_t- x_{t|t}^i\|^2.$$
\item The average root mean square error at node $i$ over the time horizon:
$$ \mathrm{RMSE}_i=\frac{1}{250}\sum_{t=1}^{250} \|x_t- x_{t|t}^i\|^2.$$ 
\item The transmission rate across the network at time $t$, i.e. the faction of nodes that transmit their data at time $t$. 
\end{itemize} 
Figure \ref{fig:rmse} and  Figure \ref{fig:rmsen} show the two aforementioned root mean square errors averaged over the runs, while Figure \ref{fig:erre} shows the corresponding averaged transmission rate. As we can see the robust filters outperform \textsf{DKF1} and \textsf{DKF2}. In particular, even in the case we increase the transmission rate, i.e. as in  \textsf{DKF2}, \textsf{RDKF} and \textsf{RDKFLOC} outperform the standard algorithm. Clearly, \textsf{DKF2} outperforms \textsf{DKF1}  in the steady state because the latter is penalized by the low transmission rate across the network. Finally,  \textsf{RDKFLOC} is slightly better than \textsf{RDKF}: it requires a lower transmission rate and it exhibits a better RMSE in the nodes characterized by a large estimation error, see Figure \ref{fig:rmsen}. This is due by the fact that the local filters of the sensor nodes in \textsf{RDKF} are too conservative, indeed their tolerances are larger than the ones in \textsf{RDKFLOC}, see Figure \ref{fig:tolloc}.

Let the average risk sensitivity parameters across the communication nodes of \textsf{RDKF} be defined as:
 \al{\theta_t^c&:=\frac{1}{|\Cc|}\sum_{i\in\Cc}\theta_t^i,\quad 
 \bar \theta_t^c:=\frac{1}{|\Cc|}\sum_{i\in\Cc}	\bar \theta_t^i\nn} 
 where $\Cc=\Nc \setminus \Sc$ and $|\Cc|$ denotes the cardinality of set $\Cc$. The average risk sensitivity parameters across the communication nodes in \textsf{RDKFLOC}, denoted by $\theta_t^{LOC,c}$ and $ \bar \theta_t^{LOC,c}$, are defined likewise.
 Figures \ref{fig:thetac} shows the aforementioned quantities averaged over the Monte Carlo runs.
 We can notice that $\theta_t^c\geq \bar \theta_t^c$ and $\theta_t^{LOC,c}\geq \bar \theta_t^{LOC,c}$, moreover we have checked that $\theta_t^i\geq \bar \theta_t^i$ for many communication nodes both in \textsf{RDKF} and \textsf{RDKFLOC}. Since the mapping $\theta\mapsto \gamma( \Omega,\theta)$ is monotone increasing, the mapping  $\Omega\mapsto \gamma( \Omega,\theta)$ is monotone decreasing according to the partial order of positive definite matrices, see   \cite{CONVTAU,LEVY_ZORZI_RISK_CONTRACTION}, and in view of the fact that
 $$   \gamma( \Omega_{t+1|t}^i,\theta_t^i)=\gamma( \bar \Omega^i_{t+1},\bar \theta_t^i),$$
 it follows that $\Omega_{t+1|t}^i\geq  \bar \Omega^i_{t+1}\geq  \bar \Psi^i_{t+1}$ for many nodes both in  \textsf{RDKF} and \textsf{RDKFLOC}. The latter inequality means that the transmission from node $i$ typically produces an increase of information in the corresponding out-neighbors nodes, which is the expected scenario. Notice that $\theta_t^c\geq \theta_t^{LOC,c}$ and $\bar \theta_t^c\geq \bar \theta_t^{LOC,c}$ which is just a consequence of the fact that \textsf{RDKF} is more conservative than \textsf{RDKFLOC}; indeed, recall that $b^i\leq b$ $\forall \,i\in\Sc$ see Figure \ref{fig:tolloc}.

 Let the average risk sensitivity parameters across the sensor nodes of \textsf{RDKF} be defined as:
 \al{\theta_t^s&:=\frac{1}{|\Sc|}\sum_{i\in\Sc}\theta_t^i,\quad 
 \bar \theta_t^s:=\frac{1}{|\Sc|}\sum_{i\in\Sc}	\bar \theta_t^i.\nn} 
 The average risk sensitivity parameters across the sensor nodes in \textsf{RDKFLOC}, denoted by $\theta_t^{LOC,s}$ and $ \bar \theta_t^{LOC,s}$, are defined likewise. 
Figure \ref{fig:thetas} shows the aforementioned quantities averaged over the Monte Carlo runs. The observations done before nodes hold also in this case. 

\section{Conclusion}\label{sec:conclusion}
In this paper we have considered the problem to estimate the state over a sensor network under model uncertainty and communication constraints. We have proposed two robust distributed strategies with event-triggered communication. More precisely, the sensor nodes compute their state estimate by solving a minimax game: one player (i.e. the estimator) aims to minimize the estimation error, while the other player selects the model in the ambiguity set  which maximizes such error. The communication among nodes is  
governed by a data-driven rule which essentially allows the data transmission   only in the case the latter provides a substantial increase of information in the nodes receiving it.  
The difference between the two distributed strategies is the way the ambiguity sets are formed at each sensor node. A stability analysis of the algorithms has been carried out showing that it is guaranteed mean-square boundedness of the state estimation error in all the nodes, under the global least favorable model, provided that the network is strongly connected, the system collectively observable and the tolerance of the ambiguity set corresponding to the global model is sufficiently small. Finally, a numerical experiment showed  that the proposed strategies are effective in the case there is model uncertainty.

\appendix
\lem\label{lemma_su_gamma} Let $\Omega>0$ and $b,\theta>0$ such that $\gamma(\Omega,\theta)=b$. Then, there exists a constant $\mu>0$ such that 
\al{\Omega-\theta I\geq\mu\Omega. \nn}
\elem
\begin{proof} The constraint $\gamma(\Omega,\theta)=b$ can be written as 
\al{\label{gamma_X}\tr(X)-\log \det(X)-n=2b}
where $X:=(I-\theta\Omega^{-1})^{-1}>0$. Let $\lambda_k>0$ denote the $k$-th eigenvalue of $X$, then (\ref{gamma_X}) can be written as 
\al{\label{gamma_X2}\sum_{k=1}^n \lambda_k -\log \lambda_k-1=2b} and the terms $ \lambda_k -\log \lambda_k-1$ in the summation are nonnegative. Hence, condition (\ref{gamma_X2}) implies that 
\al{\label{cond_f}f(\lambda_k):=\lambda_k -\log \lambda_k-1\leq 2b, \quad k=1\ldots n. } Notice that $f$ is continuous for $\lambda>0$ and it is not difficult to see that
\al{\underset{\lambda\rightarrow 0^+}{\lim}f(\lambda)=\infty,\quad \underset{\lambda\rightarrow \infty}{\lim}f(\lambda)=\infty,\quad \underset{\lambda>0}{\mathrm{argmin}} f(\lambda)=1.\nn} Accordingly, there exists $\overline\lambda >1$, which only depends on $b$, such that condition (\ref{cond_f}) is satisfied for $1\leq \lambda_k\leq \overline \lambda$, see Figure \ref{fig_graph}.
\begin{figure}[!h]\centering
\includegraphics[width=0.3\textwidth]{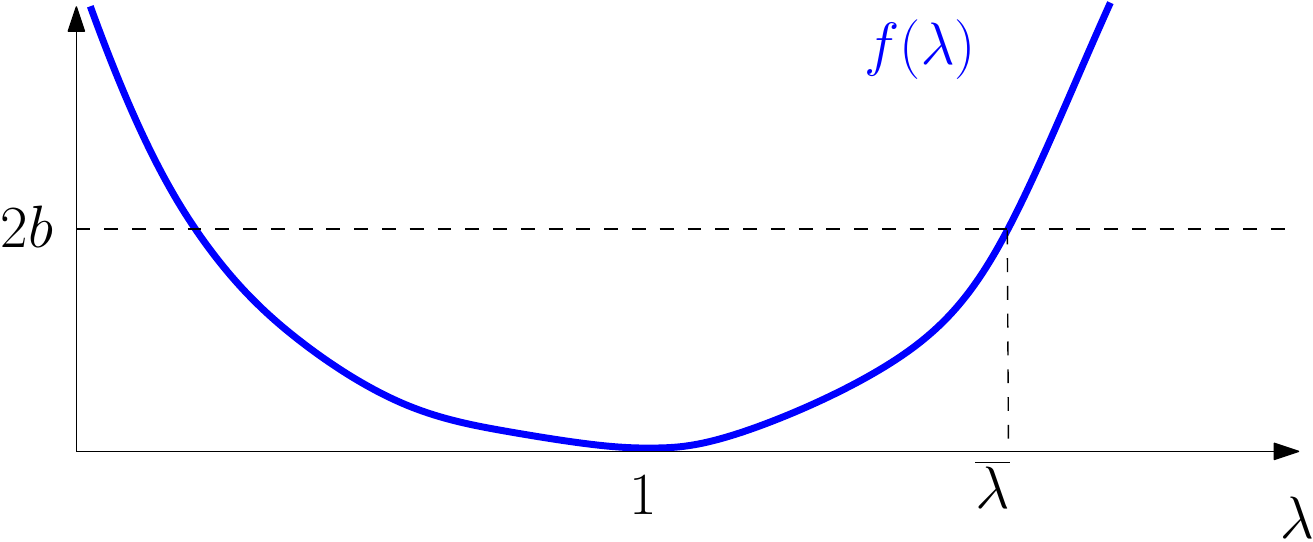} 
\caption{Pictorial description of $\overline \lambda$.}\label{fig_graph}
\end{figure} 
Accordingly, the constraint in (\ref{gamma_X}) implies that 
\al{&(I-\theta \Omega^{-1})^{-1}=X\leq \overline \lambda I\nn\\
& I-\theta \Omega^{-1}\geq \overline \lambda^{-1} I\nn\\
& \Omega-\theta I\geq \overline \lambda^{-1} \Omega\nn} and thus $\mu=\overline \lambda^{-1}$.
\end{proof}





\end{document}